\documentclass{vki}

\usepackage{lmodern,pdftexcmds,amsmath}

\usepackage{amsfonts,amsthm,bm} 

\DeclareMathAlphabet{\mathsfbi}{OT1}{\sfdefault}{bx}{sl}
\DeclareMathVersion{sfletters}
\SetSymbolFont{letters}{sfletters}{OML}{ntxsfmi}{b}{it}

\makeatletter
\newcommand{\mathbfsbilow}[1]{%
	\text{\mathversion{sfletters}$\m@th#1$}%
}

\DeclareRobustCommand{\tensor}[1]{%
	\begingroup
	\ifcat\noexpand #1\relax
	\edef\greek@test{\detokenize{#1}}%
	\edef\greek@test{\expandafter\@cdr\greek@test\@nil}%
	\edef\greek@test{\expandafter\@car\greek@test\@nil}%
	\edef\x{\the\lccode\expandafter`\greek@test}%
	\edef\y{\number\expandafter`\greek@test}%
	\ifnum\x=\y\relax
	\mathbfsbilow{#1}%
	\else
	\mathsfbi{#1}%
	\fi
	\else
	\mathsfbi{#1}%
	\fi
	\endgroup
}
\makeatother

\usepackage[most]{tcolorbox}
\usepackage{booktabs}
\usepackage{amsmath}
\usepackage{wrapfig}
\usepackage[linesnumbered,ruled,vlined]{algorithm2e}

\usepackage{mathtools}
\tcbset{enhanced,colback=cyan!2!white,colframe=cyan!90!white,coltitle=blue!20!black,fonttitle=\bfseries}

\usepackage{listings}
\usepackage{hyperref}
\hypersetup{
	colorlinks=true,
	linkcolor=blue,
	filecolor=magenta,      
	citecolor=blue,
	urlcolor=cyan,
}

\begin{document}


\pagenumbering{roman}
\title{Generalized and Multiscale Modal Analysis}
\author{Miguel Alfonso Mendez\\von Karman Institute for Fluid Dynamics 
\thanks{mendez@vki.ac.be}}

\date{12 February 2020}
\maketitle

This chapter describes modal decompositions in the framework of matrix factorizations. We highlight the differences between classic space-time decompositions and 2D discrete transforms and discuss the general architecture underpinning \emph{any} decomposition. This setting is then used to derive simple algorithms that complete \emph{any} linear decomposition from its spatial or temporal structures (bases). Discrete Fourier Transform, Proper Orthogonal Decomposition (POD), Dynamic Mode Decomposition (DMD), and Eigenfunction Expansions (EF) are formulated in this framework and compared on a simple exercise. Finally, this generalization is used to analyze the impact of spectral constraints on the classical POD, and to derive the Multiscale Proper Orthogonal Decomposition (mPOD). This decomposition combines Multiresolution Analysis (MRA) and POD. This chapter contains four exercises and two tutorial test cases. The  \textsc{Python} scripts associated to these are provided in the book's website\footnote{\url{https://www.datadrivenfluidmechanics.com/download/book/chapter8.zip}}.

\pagenumbering{arabic}
\setcounter{page}{1}
\clearpage{\pagestyle{empty}} 

\tableofcontents
\clearpage{\pagestyle{empty}}

\section{The Main Theme}\label{sec8_1}

Modal analysis aims at decomposing a dataset as a linear combination of elementary contributions called \emph{modes}. This provides the foundations to many areas of applied mathematics, including pattern recognition and machine learning, data compression, filtering, and model order reduction. A mode is the representation of the data along an element of a \emph{basis}, and the decomposition is the \emph{projection} of the dataset onto the basis (see Chapter 4). A basis maps the data onto a different space from the one in which it is sampled, with the hope that this new space is capable of better highlighting the features of interest for the task at hand. This is the central theme of this chapter: modal decompositions are ways of representing the \emph{same} information in different bases.

There is no `perfect' basis. Although some decompositions are more flexible than others, they all have their merit for a specific task. In other words, every basis can highlight or `isolate' a specific pattern in the data. Depending on how these features are used, one enters in various areas of applied mathematics-- all relevant in fluid mechanics -- as discussed in different chapters of this book.

In data compression or model order reduction, for example, one aims at representing the data with the least number of modes, to distill its `most important' features. We thus look for the smallest basis, i.e. the smallest subspace that can still handle all the information we need. 
Yet, the definition of `most important' strongly depends on the application and the experience of the user. In some settings, it is essential to identify (and retain) information contained within a specific range of frequencies. In others, it might be necessary to focus on the information that is localized within a particular location in space or time; in still others, it might be essential to extract some `coherency' or `energy' or `variance' contribution. More generally, one might be interested in a combination of the previous, and the tools described in this chapter give the reader full control over the full spectrum of options.

\subsection{The scope of modal decompositions}\label{sec8_1_1}

Data-driven modal decompositions were historically driven in fluid mechanics by the quest for identifying (and objectively define) coherence structures in turbulent flows (see Chapter 2). In parallel to this, another critical driver has always been reduced-order modeling, i.e. quest for identifying the `right' basis onto which project partial differential equations. If such a basis consists of a few carefully chosen modes, the projected system preserves the relevant dynamics and can be simulated at a fraction of the computational cost. This opens the path to fast analysis and prediction, optimization, system identification, control, and more modeling: reducing the dimensionality of a system, we facilitate regression problems such as, for example, the derivation of new models for turbulence.

In a book on data science and machine learning, this chapter could be titled \emph{linear dimensionality reduction}; this is a fundamental tool to simplify classification or regression problems. A complementary chapter should include nonlinear decompositions based on kernels or artificial neural networks (ANN), or cluster-based methods. These methods were initially developed in computer vision, data science, and statistics, and are now revolutionizing the toolbox of the fluid dynamicists, as discussed in Chapters 1, 3 and 14.

This chapter is about \emph{linear} decompositions, in which a \emph{linear combination} of modes describes the data. Each of mode might feature a nonlinear dynamic\footnote{In the sense that its evolution might not be described by a linear equation.}, but their combination is linear, with coefficients determined using inner products (see Chapter 4). In this chapter, we assume that the data is sampled in space and time, hence every mode has a spatial structure and a temporal structure. These structures represent, respectively, a basis for the discrete space and the discrete time. Furthermore, the decompositions presented consider datasets composed of a single signal. Extended methods that construct hybrid bases combining different signals (e.g. velocities and temperatures) are described in Chapter 9. 

\subsection{How many decompositions?}\label{sec_8_1_2}

Modal Decompositions can be classified into data-independent and data-driven. In the first category, the basis of the decomposition is constructed regardless of the dataset at hand and is defined by the size of the data and the basis construction criteria. This is the case of the Discrete Fourier Transform (DFT, see Chapter 4) or the Discrete Wavelet Transforms (DWT, see Chapter 5). In the first, structures are harmonics with a frequency that is an integer multiple of a fundamental one; in the second, structures are obtained by scaling and shifting a template basis (mother and father wavelets).

In data-driven decompositions, the basis is tailored to the data. The most classic examples are the Proper Orthogonal Decomposition (POD, see Chapter 6 and the Dynamic Mode Decomposition (DMD, see Chapter 7). Both POD and DMD have many variants, from which we can identify two categories of data-driven decomposition: those arising from the POD are `energy-based', those arising from the DMD are `frequency-based'. 

The POD basis is obtained from the eigenvalue decomposition of the temporal or the spatial correlation matrices. This is dictated by a constrained optimization problem that maximizes the energy (i.e., variance) along its basis, constrained to be \emph{orthogonal}, such that the error produced by an approximation of rank $\tilde{r}<R$ is the least possible. Variants of the POD can be constructed from different choices of the inner product or in the use of different averaging procedures in the computation of the correlations. Examples of the first variants are proposed by \cite{Maurel12001}, \cite{Rowley2004}, \cite{Lumley1997}, where multiple quantities are involved in the inner product. Examples of the second variants are proposed by \cite{CITRINITI2000} and \cite{towne2018spectral}, where the correlation matrix is computed in the frequency domain using time averaging over short windows, following the popular Welch's periodogram method \cite{Welch1967}. 

Within the frequency-based decompositions, the DMD basis is constructed assuming that a linear dynamical system can well approximate the data. The DMD is thus essentially a system identification procedure that aims at extracting the eigenvalues of the linear system that best fit the data. Variants of the DMD propose different methods to compute these eigenvalues. Examples are the sparsity promoting DMD \citep{Jovanovic2014pof}, the optimized DMD \citep{Chen2012}, or the randomized DMD \citep{Erichson2019}, while higher-order formulations have been proposed by \cite{Clainche2017}. Although the DMD represents the most popular formalism in fluid mechanics for such linear system identification process, analogous formulations (with slightly different algorithms) were introduced in the late '80s in climatology under the names of Principal Oscillation Patterns (POP, see \cite{Hasselmann1988} and \cite{Storch1990}) or Linear Inverse Modeling (LIM, see \cite{Penland1996}, \cite{LIM1}). 


Both `energy-based' and `frequency-based' methods have limits, illustrated in some of the proposed exercises of this chapter. These limits motivate the need for hybrid decompositions that mix the constraints of energy optimality and spectral purity. Examples of such methods for stationary datasets are the spectral proper orthogonal decomposition proposed by \cite{sieber2016spectral}, the multiresolution Dynamic Mode Decomposition \citep{MultiDMD}, the Recursive Dynamic Mode Decomposition \citep{Noack2016jfm} or the Cronos-Koopman analysis \citep{Camilleri}. A hybrid method that does not hinge on the stationary assumption is the Multiscale Proper Decomposition \citep{Mendez_Journal_2,mendez2019multi}. 

All the decompositions mentioned thus far have a common underlying architecture. This chapter presents this architecture and the formulation of the most Multiscale Proper Orthogonal Decomposition. This decomposition offer the most general formalism, unifying the energy-based and the frequency-based approaches.

\section{The General Architecture}\label{sec_8_2}

All the modal decomposition introduced in the previous section can be written as a special kind of matrix factorization. 
This view allows for defining a general algorithm for modal decomposition that is presented in Section \ref{sec_8_3}. First, Section \ref{sec_8_2_0} briefly reviews the notation followed throughout the chapter while Sections \ref{sec_8_2_1} and \ref{sec_8_2_2} put this factorization in a more general context of 2D transforms. Section \ref{sec_8_3} briefly discuss the link between discrete and continuous domain, which is essential to render all decompositions statistically convergent to a grid independent results.

\subsection{A note on notation and style}\label{sec_8_2_0}

The matrix to be factorized in this chapter is denoted as $\bm{D} (\mathbf{x_i},t_k)=\bm{D} [\mathbf{i},k]$ and is assumed to be a collection of samples of a real quantity (e.g. gray scale levels in a set of images, pressure fields in a CFD simulation or deformation fields in solid mechanics) along a spatial discretization $\mathbf{x_i}$ and a time discretization $t_k$. The notation introduced in Chapter 4 is maintained for continuous and discrete signals. The boldface notation in the spatial discretization indicates that this can be high-dimensional, e.g. $\mathbf{x_i}=(x_i,y_j)$ in a 2D domain and the boldface index $\mathbf{i}$ denotes a linear matrix index.

The matrix linear index is important when we transform a spatial realization which has the form of a matrix (for example a scalar pressure field $p[i,j]$ with $i\in[0,n_x-1]$ and $j\in[0,n_y-1]$) into a vector (that is $p[\mathbf{i}]$ with $\mathbf{i}\in[0,n_s-1]$, being $n_s=n_x\,n_y$). This index access entries of a matrix in a different way depending of whether the flattening is performed column-wise or row-wise. For example, consider the case of a matrix $\bm{A}\in\mathbb{C}^{3\times 3}$. The column-wise and the row-wise matrix indices are\footnote{Recall that here use a Python-like indexing. Hence the first entry is 0 and not 1}

$$ \mbox{column-wise }\mathbf{i}:\,\,\bm{A}=\begin{bmatrix}0 & 3 &6 \\1 & 4 &7 \\2 & 5 &8  \end{bmatrix} \quad \mbox{row-wise }\mathbf{i}:\,\,\bm{A}=\begin{bmatrix}0 & 1 &2 \\3 & 4 &5 \\6 & 7 &8  \end{bmatrix}$$

In what follows we consider a column-wise flattening. 

All the material presented in this chapter assumes a 2D space domain, the generalization to higher or lower dimensions being trivial.  Moreover, we assume that space and time domain are sampled on uniform meshes; the generalization to non-uniform mesh requires extra care in the normalization process, as discussed in Section \ref{sec_8_2_4}.
The space domain is sampled over a grid $(x_i,y_i)\in\mathbb{R}^{n_x\times n_y}$, with $x_i=i\Delta x$, $i\in[0,n_x-1]$ and $y_j=j\Delta x$, $y\in[0,n_y-1]$.

In a vector quantity, e.g. a velocity field $\bm{U}[u(\mathbf{x_i}),v(\mathbf{x_i})]$, we consider that the reshaping stacks all the components one below the other producing a state vector of size $n_s=n_C\,n_x\,n_y$ where $n_C=2$ is the number of velocity components.
Therefore, the snapshot of the data at a time $t_k$ is a vector $\mathbf{d}_k[\mathbf{i}]\in \mathbb{R}^{n_s\times 1}$ while the temporal evolution of the data at a location $\mathbf{i}$ is a vector $\mathbf{d}_{\mathbf{i}}[k]\in \mathbb{R}^{n_t\times 1}$. We then have $\bm{D}[\mathbf{i},k=c]=\mathbf{d}_c[\mathbf{i}]\in \mathbb{R}^{n_s\times 1}$ and $\bm{D}[\mathbf{i}=\mathbf{c},k]=\mathbf{d}_\mathbf{c}[k]\in \mathbb{R}^{1\times n_t}$.

\subsection{Projections in Space \emph{or} Time}\label{sec_8_2_1}

Every linear operation on a discrete 1D signal, represented as a column vector, can be carried out via matrix multiplication. This is true for convolutions, change of bases (i.e., linear transforms) and filtering. Understanding the matrix formalism of these operations brings at least three benefits. First, the notation is simplified by replacing summations on sequences with matrix-vector multiplications. Second, the efficiency and compactness of computer codes for signal processing is largely augmented\footnote{especially in interpreted languages such as Matlab or Python.}, by avoiding the nested `for' loops implied in the summation notation. Third, the linear algebra representation of these operations offers a valuable geometrical interpretation.

Consider the problem of decomposing a column vector $\bm{u}\in\mathbb{R}^{n_u\times 1}$ into a linear combination of basis vectors $\bigl\{\mathbf{b}_1,\mathbf{b}_2,...\, \mathbf{b}_{n_{b}}\bigr\}$ of equal dimensions. This means 

\begin{equation}
	\label{s8_eq1}
	u[k]=\bm{u}=\sum^{n_u}_{r=1} c_r \, \bm{b}_{r} \Longleftrightarrow \mathbf{u}=\bm{B}\,\mathbf{u}_B
\end{equation} where $\bm{B}=[\mathbf{b}_1,\mathbf{b}_2,\dots \textbf{b}_{n_{B}}]\in\mathbb{C}^{n_u\times n_b}$ is the basis matrix having all the elements of the basis along its columns and $\mathbf{u}_B=[c_1,c_1,... c_{n_b}]^T$ is the set of coefficients in the linear combination, i.e. the representation of the vector in the new basis. Computing the transform of a vector with respect to a basis $\bm{B}$ means solving a linear system of algebraic equations. Such a system can have no solution, one solution, or infinite solutions depending on $n_b$ and $n_u$.

If $n_b<n_u$, as it is the case in model order reduction, the system is \emph{over-determined} and \emph{there is no solution}\footnote{unless the $\bm{u}$ lays in the range of the column space of $\bm{B}$, in which case an unique solution exists. This is rarely the case in practice.}. In this case, we look for the approximated solution $\tilde{\mathbf{u}}_B$ that is obtained by projecting $\mathbf{u}$ onto the column space of $\bm{B}$. This is provided by the well known least-squares approximation, which gives\footnote{multiply the system by $\bm{B}^{\dag}$ and then inverts the resulting $\bm{B}^\dag \bm{B}$.}:

\begin{equation}
	\label{s8_eq2}
	\mathbf{u}=\bm{B}\,\mathbf{u}_B \Longrightarrow \tilde{\mathbf{u}}_B=(\bm{B}^\dag \bm{B})^{-1} \bm{B}^\dag \mathbf{u} \Longrightarrow
	\tilde{\mathbf{u}}=\bm{B}(\bm{B}^\dag \bm{B})^{-1} \bm{B}^\dag \mathbf{u}= \mathcal{P}_{B} \mathbf{u}\,.
\end{equation} 

The least-square solution minimizes$||\mathbf{u}-\bm{B}\,\mathbf{u}_B||_2=||\mathbf{e}||_2$; the minimization imposes that the error vector $\mathbf{e}$ is orthogonal to the column space of $\bm{B}$. In the machine learning terminology, the matrix $\mathcal{P}_{B}=\bm{B}(\bm{B}^\dag \bm{B})^{-1} \bm{B}^\dag$ is an \emph{autoencoder} which maps a signal in $\mathbb{R}^{n_u\times 1}$ to $\mathbb{R}^{n_b\times 1}$ (this is an \emph{encoding}) and then back to $\mathbb{R}^{n_u\times 1}$ (this is the \emph{decoding}).

Because the underlying linear system has no solution, the linear encoding does not generally admit an inverse if $n_b<n_u$: it is not possible to retrieve $\mathbf{u}$ from $\tilde{\mathbf{u}}_B$-- i.e. the auto-encoding loses information. A special case occurs if $n_b=n_u$. Under the assumption that the columns are linearly independent-- i.e., $\bm{B}^{-1}$ exists-- there is only one solution. It is easy to show that\footnote{use the distributive property of the inversion to show that $\mathcal{P}_{B}=\bm{B}(\bm{B}^\dag \bm{B})^{-1} \bm{B}^\dag=\bm{B}\bm{B}^{-1} \bigl (\bm{B}^\dag\bigr)^{-1}\bm{B}=\bm{I}$.} this yields $\mathcal{P}_{B}=\mathbf{I}$: in this case $\tilde{\mathbf{u}}=\mathbf{u}$, and the tilde is removed because the auto-encoding is lossless. This is fundamental in filtering applications, for which the basis matrix is usually square: a signal is first projected onto a certain basis (e.g., Fourier or Wavelets), manipulated, and then projected back. 

The last possibility, that is $n_b>n_u$, results in an \emph{under-determined} system. In this case, there are \emph{infinite solutions}. Among these, it is common practice to consider the one that leads to the least energy in the projected domain, i.e., such that $min(||\mathbf{u}_B||)$. This approach, which also yields a reversible projection, is known as the least norm solution and reads:
\begin{equation}
	\label{s8_eq3}
	\mathbf{u}=\bm{B}\,\mathbf{u}_B \Longleftrightarrow {\mathbf{u}}_B=\bm{B}^{\dag}(\bm{B} \bm{B}^\dag)^{-1}  \mathbf{u} \,.
\end{equation}

It is now interesting to apply these notions to decompositions (projections) in space and time domains.
Considering first the projection in the space domain, let the signal in \eqref{s8_eq2} be a column of the dataset matrix, i.e. $\mathbf{u}:=\mathbf{d}_k[\mathbf{i}]\in\mathbb{R}^{n_s\times 1}$. Let $\bm{\Phi}=[\bm{\phi}_1,\bm{\phi}_2,\dots, \bm{\phi}_{n_{\phi}}] \in \mathbb{C}^{n_{s}\times n_{\phi}}$ denote the spatial basis. Because the matrix multiplication acts independently on the columns of  $\bm{D}$, equation \eqref{s8_eq1} is

\begin{equation}
	\label{s8_eq4}
	\mathbf{d}_k=\bm{\Phi}\,\mathbf{d}_{\phi} \Rightarrow \tilde{\mathbf{d}}_{k \phi}=(\bm{\Phi}^\dag \bm{\Phi})^{-1} \bm{\Phi}^\dag \mathbf{d}_k \,\,\,\mbox{;}\,\,\, 
	\bm{D}=\bm{\Phi}\bm{D}_{\phi}\,\Rightarrow \tilde{\bm{D}}_{\phi}=(\bm{\Phi}^\dag \bm{\Phi})^{-1} \bm{\Phi}^\dag \bm{D}.
\end{equation}

Here the transformed vector is $\tilde{\mathbf{d}}_{k \phi}$ while the matrix $\tilde{\bm{D}}_{\phi}$ collects all the transformed vectors, i.e. the coefficients of the linear combinations of basis elements $\{\bm{\phi}_1, \bm{\phi_2},\dots ,\bm{\phi}_R\}$ that represents a given snapshot $\bm{d}_k$. 

The same reasoning holds for transforms in the time domain. In this case, let the signal be a row of the dataset matrix, i.e. $\mathbf{u}:= \mathbf{d}^T_{\mathbf{i}}[k]\in\mathbb{R}^{n_t\times 1}$. Defining the temporal basis matrix as $\bm{\Psi}=[\bm{\psi}_1,\bm{\psi}_2,\dots, \bm{\psi}_{n_{\psi}}] \in \mathbb{C}^{n_t\times n_{\psi}}$ and handling the transpositions with care\footnote{More generally, in a matrix multiplication $\bm{A}\bm{B}$, the matrix $\bm{A}$ is acting on the columns of $\bm{B}$. The same action along the rows of $\bm{B}$ is obtained by $\bm{A}\bm{B}^T$.}, the analogous of \eqref{s8_eq4} in the time domain reads:

\begin{equation}
	\label{s8_eq5}
	\mathbf{d}^T_{\mathbf{i}}=\bm{\Psi}\,\mathbf{d}^T_{\psi} \Rightarrow \tilde{\mathbf{d}}^T_{\mathbf{i} \psi}=(\bm{\Psi}^\dag \bm{\Psi})^{-1} \bm{\Psi}^\dag \mathbf{d}^T_{\mathbf{i}} \,\,\,\mbox{;}\,\,\,  
	\bm{D}=\bm{D}_{\psi}\bm{\Psi}^T\,\Rightarrow \tilde{\bm{D}}_{\psi}=\bm{D}\overline{\bm{\Psi}}(\bm{\Psi}^\dag {\bm{\Psi}})^{-1}\,\,.
\end{equation}

Here the transformed vectors are $\tilde{\mathbf{d}}_{\mathbf{i} \psi}$ and the matrix $\tilde{\bm{D}}_{\psi}$ collects in its columns all the transformed vectors, i.e. the coefficients of the linear combinations of basis elements $\{\bm{\psi}_1, \bm{\psi_2},\dots ,\bm{\psi}_R\}$ that represents a the data evolution at a location $\mathbf{d}_{\mathbf{i}}[k]$.

Note that approximations can be obtained along space and time domains as

\begin{equation}
	\label{s8_eq6}
	\tilde{\bm{D}}_{\phi}=\bm{\Phi}(\bm{\Phi}^\dag \bm{\Phi})^{-1} \bm{\Phi}^\dag \bm{D}\quad \mbox{and} \quad \tilde{\bm{D}}_{\psi}=\bm{D} \overline{\bm{\Psi}}(\bm{\Psi}^\dag \bm{\Psi})^{-1} \bm{\Psi}^T\,
\end{equation}

These approximations are exact (the projections are reversible) if the bases are complete, i.e. the matrices $\bm{\Phi}$ and $\bm{\Psi}$ are square\footnote{and, of course, have full rank.} ($n_{\phi}=n_s$ and $n_{\psi}=n_t$). In these cases, autoencoding is loseless: $\tilde{\bm{D}}_{\phi}=\bm{D}$ and $\tilde{\bm{D}}_{\psi}=\bm{D}$.

It is left as an exercise to see how \eqref{s8_eq4}-\eqref{s8_eq5}-\eqref{s8_eq6} simplify if the bases in the space and time are also orthonormal, i.e. the inner products yields $\bm{\Phi}^\dag \bm{\Phi}=\mathbf{I}$ and $\bm{\Psi}^\dag \bm{\Psi}=\mathbf{I}$ regardless of the number of basis elements\footnote{If $\bm{\Psi}^\dag \bm{\Psi}=\mathbf{I}$, then taking the conjugation on both sides gives $\bm{\Psi}^T \overline{\bm{\Psi}}=\mathbf{I}$ while taking the transposition on both sides gives $\overline{\bm{\Psi}} \bm{\Psi}^T=\mathbf{I}$.} $n_{\phi}$ and $n_{\psi}$.  On the other hand, from the fact that \eqref{s8_eq6} is \emph{exact} only for \emph{complete} bases, one can see that $\bm{\Phi}\bm{\Phi}^\dag=\mathbf{I}$ \emph{only if} $n_{\phi}=n_s$ and $\bm{\Psi} \bm{\Psi}^\dag=\mathbf{I}$ \emph{only if } $n_{\psi}=n_t$.

\subsection{Projections in Space \emph{and} Time}\label{sec_8_2_2}

Let us now consider 2D transforms of the dataset matrix. These are the discrete version of 2D transforms, which in this case act on the columns and the rows of the matrix independently. Let $\bm{D}_{\phi \psi}$ be the 2D transform of the matrix $\bm{D}$ taking $\bm{\Phi}=[\bm{\phi}_1,\bm{\phi}_2,\dots, \bm{\phi}_{n_{\phi}}] \in \mathbb{C}^{n_{s}\times n_{\phi}}$ and $\bm{\Psi}=[\bm{\psi}_1,\bm{\psi}_2,\dots, \bm{\psi}_{n_{\psi}}] \in \mathbb{C}^{n_t\times n_{\psi}}$ as basis matrices for the columns and the rows respectively. The general 2D transform pair is:

\begin{equation}
	\label{s8_eq7}
	\bm{D}=\bm{\Phi} \bm{D}_{\phi \psi} \bm{\Psi}^T \Leftrightarrow \tilde{\bm{D}}_{\phi \psi}=(\bm{\Phi}^\dag \bm{\Phi})^{-1}\bm{\Phi}^\dag\,\bm{D}\,\overline{\bm{\Psi}}(\bm{\Psi}^\dag \bm{\Psi})^{-1}\,.
\end{equation} 

A classic example of this kind of transforms is the 2D Fourier transform, in which case the bases matrices are

\begin{equation}
	\label{s8_eq8}
	\bm{\Phi}_{\mathcal{F}}[m,n]=\frac{1}{\sqrt{n_s}} e^{{2 \pi \mathrm {j}\, m\,n}/{n_s}}\,\,\mbox{and} \,\,
	\bm{\Psi}_{\mathcal{F}}[m,n]=\frac{1}{\sqrt{n_t}} e^{{2 \pi \mathrm {j}\, m\,n}/{n_t}}\,.
\end{equation} 

These are the \emph{Fourier Matrices}, which are orthonormal. If these are also square (i.e., the number of frequency bins equals $n_s$ and $n_t$ respectively), these matrices are also symmetric. Then, \eqref{s8_eq7} is reversible and the 2D Discrete Fourier pair is\footnote{Observe that transpositions are removed because of the symmetry of the bases.}

\begin{equation}
	\label{s8_eq9}
	\bm{D}=\bm{\Phi}_{\mathcal{F}} \widehat{\bm{D}} \bm{\Psi}_{\mathcal{F}} \Longleftrightarrow \widehat{\bm{D}}=\overline{\bm{\Phi}}_{\mathcal{F}}\,\bm{D}\,\overline{\bm{\Psi}}_{\mathcal{F}}\,,
\end{equation}

The (discrete) Fourier transform plays an essential role in Section \ref{sec_8_4}. Let us use a wide hat $\,\widehat{}\,$ to indicate the Fourier transform of a matrix, whether this is carried out along its rows, columns, or both. 

Different decompositions can be derived by taking different bases for the rows and for the columns. This is often done in the 2D Discrete Wavelet transform (see Chapter 5), as one might desire to highlight different features along different directions. For example, in the literature of image processing, we define \emph{horizontal details} as those features that have high frequencies over the column space and low frequencies over the raw space. We can thus construct a basis (i.e design a filter) that enhance or reduce these details.

In the analysis of 2D datasets such as images, decompositions of these form common and useful. For reduced-order modeling, however, these have a substantial disadvantage: the number of basis elements (i.e., its `modes') are $n_s\times n_t$. While most of the data have a sparse representation in a 2D Fourier or Wavelet bases\footnote{i.e. most of the entries in  $\bm{D}_{\phi }$, $\bm{D}_{ \psi}$ or $\bm{D}_{\phi \psi}$ are zero or almost zero}, this representation is still inefficient. The inefficiency stems from the lack of a `variable separation': the entry $\bm{D}_{\phi \psi}[m,n]$ measure the correlation with a basis matrix constructed from the $m-th$ basis element of the column space and the $n-th$ basis element of the row space. Even in the ideal case of orthonormal bases on both columns and  rows (i.e., space and time), the decomposition in \eqref{s8_eq7} requires a double summation:

\begin{equation}
	\label{s8_eq10}
	{\bm{D}}[\mathbf{i},k]=\sum^{n_s-1}_{m=0}\sum^{n_t-1}_{n=0}\,\bm{D}_{\phi \psi}[m,n]\bm{\phi}_m[\mathbf{i}]\,\bm{\psi}_n[k]
\end{equation}

In modal analysis, we seek a transformation that renders $\bm{D}_{\phi \psi}$ \emph{diagonal} and that has no more than $R=min(n_s,n_t)$ modes. We seek separation of variables and hence enter data-driven modal analysis from the next section.

\subsection{The Fundamental Factorization}\label{sec_8_2_3}

We impose that $\bm{D}_{\phi \psi}$ is diagonal. In this case, the transformed matrix is denoted as $\bm{\Sigma}=diag(\sigma_1,\dots,\sigma_R)\in \mathbb{R}^{R\times R}$ and the entries $\sigma_r$ are referred to as \emph{mode amplitudes}. The summation in \eqref{s8_eq10} and the factorization in \eqref{s8_eq7} become

\begin{equation}
	\label{s8_eq11}
	\bm{D}=\bm{\Phi} \,\bm{\Sigma}\,\bm{\Psi}^T \Longleftrightarrow \bm{D}[\mathbf{i},k]=\sum^{R-1}_{r=0} \sigma_r \bm{\phi}_m[\mathbf{i}]\,\bm{\psi}_n[k]\Longleftrightarrow \bm{D}=\sum^{R-1}_{r=0} \sigma_r \bm{\phi}_r \,\bm{\psi}_r^T\,.
\end{equation} Each of the terms $\sigma_r \bm{\phi}_r \,\bm{\psi}_r^T$ produces a matrix of unitary rank, which we call \emph{mode}. Note that the rank of a (full rank) rectangular matrix is $rank(\bm{D})=min(n_s,n_t)$, but in general the number of relevant modal contributions are not associated to the rank\footnote{For example, we might add five modes and obtain a matrix of rank four. The number of modes equals the rank of the matrix only if these modes are orthogonal in both space and time, i.e. only if these are POD modes.}. We obtain approximations of the dataset by zeroing some of the entries along the diagonal of $\bm{\Sigma}$. As in the previous section we use tildes to denote approximations. Moreover, since infinite decompositions could be obtained by dividing the diagonal entries $\sigma_r=\bm{\Sigma}[r,r]$ by the length of the corresponding basis matrices, we here assume that $||\bm{\phi}_r||=||\bm{\psi}_r||=1$ $\forall r\in[0,\dots,R-1]$.

The reader has undoubtedly recognized that the factorization in \eqref{s8_eq11} has the same structure of the Singular Value Decomposition  (SVD) introduced in Chapter 6. The SVD is a very special case of \eqref{s8_eq11}, but there are \emph{infinite} other possibilities, depending on the choice of the basis. Nevertheless, such a choice has an important constraint that must be discussed here. 

Assume that the spatial structures $\bm{\Phi}$ are given. Projecting on the left gives:

\begin{equation}
	\label{s8_eq12}
	(\bm{\Phi}^\dag \bm{\Phi})^{-1} \bm{\Phi}^{\dag}\bm{D}=\bm{\Sigma}\,\bm{\Psi}^T=\bm{D}_{\phi}=\begin{bmatrix}
		\sigma_1 \psi_1[1]   &\sigma_1 \psi_1[2]& \dots &\sigma_1 \psi_1[n_t] \\
		\sigma_2 \psi_2[1]   &\sigma_2 \psi_2[2]& \dots &\sigma_2 \psi_2[n_t] \\
		\vdots   &\vdots&\vdots&\vdots\\
		\sigma_R \psi_R[1]   &\sigma_R \psi_R[2]& \dots &\sigma_R \psi_R[n_t]  \\
	\end{bmatrix}\,.  \end{equation}

This is equivalent to \eqref{s8_eq4}: this is a transform along the columns of $\bm{D}$. The transformed snapshot, i.e. the set of coefficients in the linear combination of $\bm{\phi}$'s, change from snapshot to snapshot -- they evolve in time. The key difference with respect to a 2D transform is that the evolution of the structure in the basis element $\bm{\phi}_r$ only depends on the corresponding temporal structure $\bm{\psi}_r$. Moreover, notice that $\sigma_r=||\sigma_r\,\psi_r||$, since $||\psi_r||=1$: this explains  while the amplitudes are real quantities by construction: they represents the norm of the $r-th$ column $\bm{D}_{\phi}$, which we denote as $\bm{D}_{\phi}[:,r]$ following a \textsc{Python} notation.

The same observation holds if the temporal structures $\bm{\Psi}$ are given. In this case, projecting on the right gives:

\begin{equation}
	\label{s8_eq13}
	\bm{D}\overline{\bm{\Psi}}(\bm{\Psi}^\dag {\bm{\Psi}})^{-1}=\bm{\Phi}\bm{\Sigma}=\bm{D}_{\psi}=\begin{bmatrix}
		\sigma_1 \phi_1[1]   &\sigma_2 \phi_1[1]& \dots &\sigma_R \phi_R[1] \\
		\sigma_1 \phi_2[2]   &\sigma_2 \phi_2[2]& \dots &\sigma_R \phi_R[2] \\
		\vdots   &\vdots&\vdots&\vdots\\
		\sigma_1 \phi_R[n_s]   &\sigma_2 \phi_R[n_s]& \dots &\sigma_R \phi_R[n_s]  \\
	\end{bmatrix}\,.  \end{equation}

This is equivalent to \eqref{s8_eq5}: this is a transform along the rows of $\bm{D}$.
The set of coefficients in the linear combination of $\bm{\psi}$'s are spatially distributed according to the corresponding spatial structures $\bm{\phi}$'s. As before, we see that the amplitudes of the modes can also be computed as $\sigma_r=||\sigma_R \bm{\phi}_r||=||\bm{D}_{\psi}[:,r]||$.

\begin{figure}[h]
	\centering
	\includegraphics[keepaspectratio=true,width=1 \columnwidth]{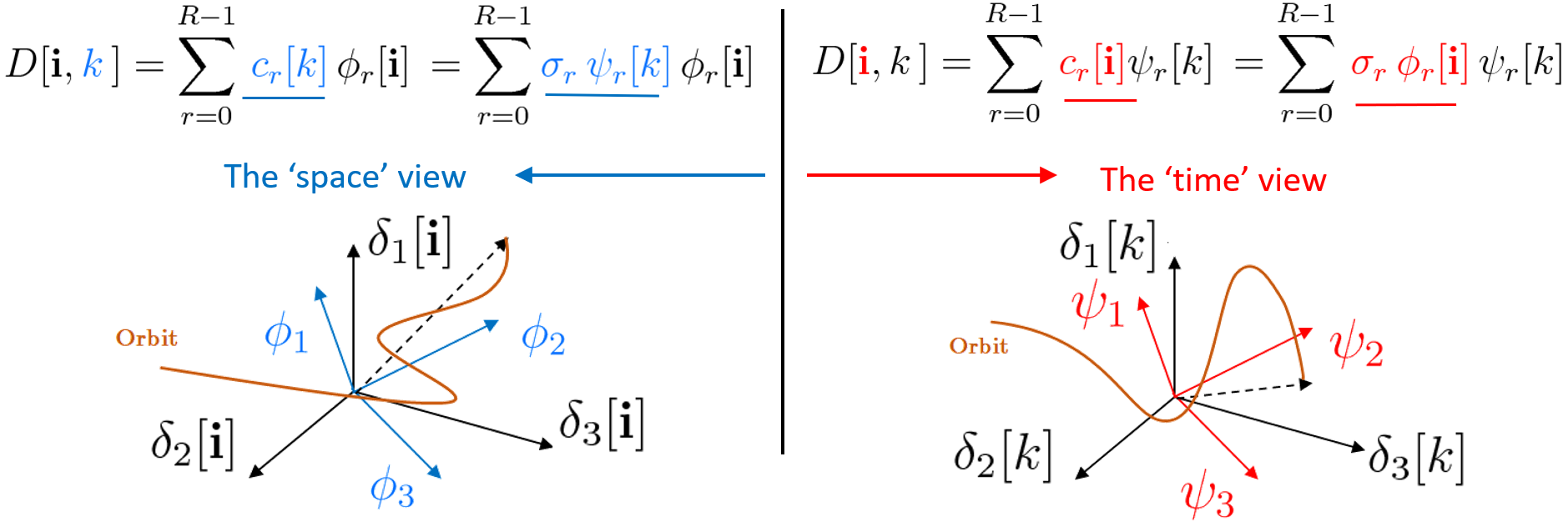}
	\caption{The `space view' and the `time-view' of modal analysis. In the `space view', every snapshot $\mathbf{d}_k[\mathbf{i}]$ is a linear combination of basis elements $\bm{\phi}$'s. The coefficients of this combination evolve in time according to the associated $\bm{\psi}$'s. In the `time-view', every temporal evolution $\mathbf{d}_{\mathbf{i}}[k]$ is a linear combination of basis elements $\bm{\psi}$'s. The coefficients of this combination are spatially distributed according to the associated $\bm{\phi}$'s.}
	\label{ch8Fig1}
\end{figure}

The space-time symmetry is further elucidated in Figure \ref{ch8Fig1}, considering also the cases in which $\bm{\Phi}:= \mathbf{I}/\sqrt{n_s}$ or $\bm{\Psi}:= \mathbf{I}/\sqrt{n_t}$. In the `space view' on the left, we follow the data in time from a specific spatial basis. If this basis is the set of impulses $\delta_k[\mathbf{i}]$, the temporal evolutions along each element of the basis is given by the time evolution  $\mathbf{d}_{\mathbf{i}}[k]$. 
In the `time view' on the right, we analyze the spatial distribution from a specific temporal basis. If this basis is the set of impulses $\delta_{\mathbf{i}}[k]$, the spatial distribution of each member of the basis (i.e., every instant) is given by the snapshots $\mathbf{d}_{k}[\mathbf{i}]$ itself.

With these views in mind, the reader should understand the most important observations of this section: \emph{in the decomposition in }\eqref{s8_eq11} \emph{it is not possible to impose both basis matrices} $\bm{\Phi}$ \emph{and} $\bm{\Psi}$-- \emph{given one of the two, the other is univocally determined}. This observation also leads to the formulations of two general algorithms to compute this factorization given one of the two bases. These are listed below:

\begin{algorithm}
	\KwIn{$\bm{D}\in \mathbb{R}^{n_s\times n_t}$ and $\bm{\Phi}\in \mathbb{R}^{n_s\times R}$}
	\KwOut{$\bm{\Sigma}\in\mathbb{R}^{R\times R}$ and $\bm{\Psi}\in\mathbb{R}^{n_t\times R}$}
	Project $\tilde{\bm{D}}_{\phi}=(\bm{\Phi}^\dag \bm{\Phi})^{-1} \bm{\Phi}^\dag\bm{D}$ (eq. 12)\\
	From \eqref{eq12} it is $\bm{D}_{\phi}=\bm{\Sigma} \bm{\Psi}^T$, hence :\\
	\text{for $r<=R$}, r=r++, $\sigma_r=||\bm{D}_{\phi}[r,:]||_2$;  $\psi_r=\bm{D}_{\phi}[r,:]/\sigma_r$ \\ 
	Assembly $\bm{\Sigma}=diag(\sigma_0,\sigma_1,\dots \sigma_R)\in \mathbb{R}^{R\times R}$ and $\bm{\Psi}=[\psi_0,\psi_1,\dots \psi_R]\in \mathbb{C}^{n_t \times R}$
	\caption{Projection from $\bm{\Phi}$}
	\label{algo:a}
\end{algorithm}

\vspace{-0.5mm}
\begin{algorithm}
	\KwIn{$\bm{D}\in \mathbb{R}^{n_s\times n_t}$ and $\bm{\Psi}\in \mathbb{R}^{n_t\times R}$}
	\KwOut{$\bm{\Sigma}\in\mathbb{R}^{R\times R}$ and $\bm{\Phi}\in\mathbb{R}^{n_s\times R}$}
	Project $\tilde{\bm{D}}_{\psi}=\bm{D}\overline{\bm{\Psi}}(\bm{\Psi}^\dag {\bm{\Psi}})^{-1}=\bm{\Phi}\bm{\Sigma}$ (eq. 13)\\
	From \eqref{eq12} it is $\bm{D}_{\psi}=\bm{\Phi}\bm{\Sigma}$, hence :\\
	\text{for $r<=R$}, r=r++, $\sigma_r=||\bm{D}_{\psi}[r,:]||_2$;  $\phi_r=\bm{D}_{\phi}[r,:]/\sigma_r$ \\ 
	Assembly $\bm{\Sigma}=diag(\sigma_0,\sigma_1,\dots \sigma_R)\in \mathbb{R}^{R\times R}$ and $\bm{\Phi}=[\phi_0,\phi_1,\dots \phi_R]\in \mathbb{C}^{n_s \times R}$
	\caption{Projection from $\bm{\Psi}$}
	\label{algo:b}
\end{algorithm}

Finally, note that these algorithms are the most general ones to complete a decomposition from its spatial or its temporal structures. Such a level of generality highlights the common structure but also leads to the least efficient approach: every decomposition offers valuable short-cuts that we discuss in Section \ref{sec_8_3}.

\subsection{Amplitudes and Energies}\label{sec_8_2_4}

The normalization of a data-driven decomposition cannot be performed beforehand, because the basis is initially unknown. We here analyze the consequences of imposing a unitary norm for the bases vectors $\bm{\phi}$ and $\bm{\psi}$ and the link of this operation with the continuous domain. Using the continuous inner product (see Chapter 4), the energies of continuous basis element in space and time are:

\begin{equation}
	\label{s8_eq14}
	||\phi_r(\mathbf{x})||^2=\frac{1}{\Omega}\int_{\Omega} \phi_r(\mathbf{x})\,\overline{\phi}_r(\mathbf{x})\, d\Omega \,\,\,\,\mbox{and} \,\,\,\, ||\psi_r(t)||^2=\frac{1}{T}\int_{T} \psi_r(t) \overline{\psi}_r(t)\, dT\,
\end{equation}\, having taken the spacial basis in a domain\footnote{this could be an area in a 2D domain or a volume in a 3D domain. In this sense, the bases in the continuous domain are always considered in terms of `densities'.} $\Omega$, and the temporal basis in a domain $T$. If a Cartesian and uniform mesh is used in space and time, and if a right endpoint integration scheme is considered, these are approximated by their discrete equivalent 
\begin{equation}
	\label{s8_eq15}
	||\phi_r(\mathbf{x})||^2=\frac{1}{n_s} ||\bm{\phi}_r||^2=\frac{1}{n_s} \bm{\phi}^{\dag}_r\bm{\phi} \,\,\,\,\mbox{and} \,\,\,\, ||\psi_r(t)||^2=\frac{1}{n_t} ||\bm{\psi}_r||^2=\frac{1}{n_t} \bm{\psi}^{\dag}_r\bm{\psi}\,,
\end{equation} where the norms are simply Euclidean norms $||\bm{a}||=\bm{a}^{\dag} \bm{a}$. For nonuniform meshes or more sophisticated integration schemes, weighted inner products must be introduced (see Chapter 6). 

Consider now an approximation of the data using only one mode. This is a matrix $\tilde{\bm{D}}[\mathbf{i},k]=\sigma_r \phi_{r}[\mathbf{i}] \psi_{r}[k]$ of unitary rank. The \emph{total} energy associated to this mode, assuming that this is the discrete version of a continuous space-time evolution, is 

\begin{equation}
	\label{s8_eq16}
	\mathcal{E}\{\tilde{D}[\mathbf{x},t]\}=\frac{1}{\Omega\, T}\int_{T}\int_{\Omega} \tilde{D}^2(\mathbf{x},t) d\Omega dt \approx \frac{1}{n_t n_s} \sum^{n_s}_{\mathbf{i}=0} \sum^{n_t}_{k=0} \tilde{\bm{D}}^2[\mathbf{i},k]\,.
\end{equation}

The summation on the right is the Frobenious norm\footnote{Recall: the Frobenious norm of a matrix $\bm{A}\in\mathbb{C}^{n_s\times n_t}$ can be written as $||\bm{A}||_F=\mbox{tr}\{\bm{A}^{\dag}\bm{A}\}=\mbox{tr}\{\bm{A}\bm{A}^{\dag}\}$.}, hence we have

\begin{equation}
	\label{s8_eq17}
	\mathcal{E}\{\tilde{D}[\mathbf{x},t]\}\approx \frac{1}{n_t n_s} \mbox{tr}\{\tilde{\bm{D}}^\dag\tilde{\bm{D}}\}=\frac{\sigma^2_r}{n_t n_s} \mbox{tr}\{ \bm{\psi}_r \bm{\phi}^\dag_r \bm{\phi}_r  \bm{\psi}^\dag_r\}=\frac{\sigma^2_r}{n_t n_s} .
\end{equation} Observe that $\bm{\phi}^\dag_r \bm{\phi}_r=1$ because of the unitary length constraint on the spatial structures and $\mbox{tr}\{\bm{\psi}_r \bm{\psi}_r^T\}=\bm{\psi}^{\dag}_r \bm{\psi}_r=1$ for the same constraint on the temporal structures. \emph{We thus conclude that the square of the normalized amplitudes $\hat{\sigma}_r=\sigma_r/\sqrt{n_t n_s}$ are mesh independent approximation of the energy contribution of each mode}.

This statement is extremely important but requires great care: \emph{in general, summing the energy contribution of each mode does not equal the energy in the entire data}. This occurs only if \emph{both the temporal and the spatial structures are orthonormal}- that is only for the Proper Orthogonal Decomposition (see Chapter 6). For any other decomposition, the energy of an approximation with $\tilde{r}<R$ modes is

\begin{equation}
	\label{s8_eq18}
	\mathcal{E}\{\tilde{D}[\mathbf{x},t]\}=\frac{1}{n_s n_t} \mbox{tr}\Bigl\{\tilde{\bm{\Psi}}\,\tilde{\bm{\Sigma}}\,\tilde{ \bm{\Phi}}^{\dag}\,\tilde{\bm{\Phi}}\,\tilde{\bm{\Sigma}}\,\tilde{\bm{\Psi}}^\dag\Bigr\}=\frac{1}{n_s n_t}
	\mbox{tr}\Bigl\{\tilde{\bm{\Phi}}\tilde{\bm{\Sigma}}\tilde{\bm{\Psi}}^\dag\,\tilde{\bm{\Psi}}\,\tilde{\bm{\Sigma}}\,\tilde{\bm{\Phi}}^{\dag}\Bigr\}
\end{equation}

The convergence of a decomposition indicates how efficiently this can represent the data set, that is how quickly the norm of the approximation converges to the norm of the data as $\tilde{r}\rightarrow R$. Convergence can be either measured in terms of $l^2$ norm or in terms of Frobenius norm. Using the $l^2$ norm, the convergence can be measured as 

\begin{equation}
	\label{s8_eq19}
	E(\tilde{r})=\frac{||\bm{D}-\tilde{\bm{D}}(\tilde{r})||}{||\bm{D}||}=\frac{||\bm{D}-\tilde{\bm{\Phi}}\tilde{\bm{\Sigma}}\tilde{\bm{\Psi}}||_2}{||\bm{D}||}
\end{equation}

\section{Common Decompositions}\label{sec_8_3}

We here describe how common decompositions fits in the factorization framework. We adopt a `time-view', i.e. we identify every decomposition from its temporal structure $\bm{\Psi}$ and complete the factorization following algorithm 2 from the previous section.

\subsection{The Delta Decomposition}\label{sec_8_3_1}

This decomposition is useless, but conceptually important. Take basis of impulses $\bm{\Psi}_{\delta}=\mathbf{I}/\sqrt{n_t}$ as the temporal basis. This basis is orthonormal and we distinguish this decomposition from the others using the subscript $\delta$.

This decomposition generalizes the convolution representation of a signal introduced in Chapter 4. The general temporal structure is $\psi_r[k]=\bm{\delta}_r[k]=\delta[k-r]$. The decomposition of the snapshot $\bm{d}_k[\mathbf{i}]$ becomes:

\begin{equation}
	\label{s8_eq20}
	\bm{d}_k[\mathbf{i}]=\sum^{R-1}_{r=0} \sigma_{\delta\,r} \bm{\phi}_{\delta r}[\mathbf{i}] \bm{\psi}_{\delta r}[k] = \sum^{R-1}_{r=0} \sigma_{\delta\,r} \bm{\phi}_{\delta r}[\mathbf{i}] \frac{\delta[k-r]}{\sqrt{n_t}}= \frac{\sigma_{\delta\,k} \bm{\phi}_{\delta k}[\mathbf{i}] }{\sqrt{n_t}}\,.
\end{equation}

Therefore the spatial basis corresponds to the normalized snapshots:

\begin{equation}
	\label{s8_eq21}
	\bm{\phi}_{\delta k}[\mathbf{i}]=\frac{\bm{d}_k[\mathbf{i}]}{||\bm{d}_k[\mathbf{i}]||} \quad \mbox{and} \quad \sigma_{\delta\,k}=||\bm{d}_k[\mathbf{i}]||\sqrt{n_t}=||\bm{d}_k(\mathbf{x_i})||\sqrt{n_t\,n_s}
\end{equation}

This decomposition yields the highest possible temporal localization: large amplitudes $\sigma_{\delta r}$ corresponds to a snapshot with large energy content. A strong amplitude decay indicates that some snapshots have much lower energy content than others. In a dataset in which a strong impulsive event occurs, the POD temporal structure may approach a delta-like form: when this happens, the associated spatial structure tends to the corresponding normalized snapshot.
A similar limit is obtained in a DWT (see Chapter 5) if the finest scale is pushed to its limit. In a Haar Wavelet transform, for example, the finest scale correspond to elements of the form $\bm{w}_r=(\bm{\delta}_r+\bm{\delta}_{r+1})/\sqrt{2}$. The spatial structures associated to these are normalized averages of the snapshots $\textbf{d}_{k}$ and $\textbf{d}_{k+1}$. Similarly, the spatial structure associated with the basis of the form $\bm{\psi}_r=1/\sqrt{n_t}$ is the normalized time average of the data.

\subsection{The Discrete Fourier Transform (DFT)}\label{sec_8_3_2}

We here refer to the DFT of the dataset in time, i.e. this is equivalent to performing the DFT of $\bm{D}$ row-wise as in \eqref{s8_eq5}, followed by a normalization to compute the amplitudes. In the DFT we set $\bm{\Psi}=\bm{\Psi}_{\mathcal{F}}$ 
and we use the subscript $\mathcal{F}$ to distinguish it from the others. This matrix is defined in \eqref{s8_eq8}.  This basis is orthonormal and symmetric. Using \eqref{s8_eq13} and the general algorithm 2 gives

\begin{equation}
	\label{s8_eq22}
	\bm{D}\overline{\bm{\Psi}}_{\mathcal{F}}=\bm{\Phi}_{\mathcal{F}}\bm{\Sigma}_{\mathcal{F}}=\widehat{\bm{D}}\Longleftrightarrow \sigma_{\mathcal{F}r}=||\bm{\phi}_r \sigma_r||=||\bm{D}\overline{\bm{\psi}}_r||\,,\quad \mbox{i.e.}\,\, \bm{\Phi}_{\mathcal{F}}= \widehat{\bm{D}} \bm{\Sigma}^{-1}_{\mathcal{F}}\,.
\end{equation} 

The first step on the left is the row-wise Fourier Transform of $\bm{D}$, denoted as $\widehat{\bm{D}}$. The second step is the normalization, which can be computed along the rows of $\widehat{\bm{D}}$, and in the last step we see that the spatial structures correspond to the normalized distribution of a given entry of the Fourier spectra. We thus see that each entry of the spatial basis tells how much the data at a specific location correlates with the harmonic associated with the discrete frequency $r$. The computation of the digital frequencies from the sampling frequency of the data is described in Chapter 4. In the same section, it was recalled that the DFT implicitly assumes that the data is periodic of period $T=n_t\Delta t $. 

From a computational point of view, the multiplication by $\bm{\Psi}_{\mathcal{F}}$ can be efficiently carried out using the FFT algorithm. Note that \textsc{fft} routines from \textsc{Matlab} or \textsc{Python} define the Fourier basis matrix as the conjugate of \eqref{s8_eq8}. This decomposition is, by definition, the one with highest frequency resolution. As described in Chapter 6, this comes at the price of no time localization, since this temporal basis is of infinite duration. An important property of the Fourier basis is that its elements are eigenvectors of circulant matrices \citep{Gray2005}. 

This property links the DFT and the POD in stationary data, and it is worth a brief discussion. An intuitive derivation follows from the notions introduced in Chapter 4, specifically the convolution theorem in the discrete domain and the link between impulse response of an LTI system and frequency transfer functions.

Hence, let $\mathbf{u}:= u[k]\in\mathbb{R}^{n_t\times1}$ with $k\in[0,n_t-1]$ be a discrete signal, represented by a column vector of appropriate size, which is introduced as an input to a SISO system. Let $\mathbf{y}:= y[k]\in\mathbb{R}^{n_t\times1}$ be the corresponding output. Let $h[k]$ be the \emph{impulse response} of the LTI system, assuming that this is of finite duration and that it has at most $n_t$ non-zero entries. The input-output link in the time domain is provided by the convolution sum (see Chapter 4). If both input and output are periodic, the convolution can be written using a Circulant matric $\bm{C}_h$ and a matrix multiplication:

\begin{equation}
	\label{s8_eq23}
	\begin{split}
		y[k]=\sum^{l=\infty}_{l=-\infty} u[k]h[k-l]\rightarrow \mathbf{y}=\bm{C}_h \mathbf{u}\quad \quad\\
		\begin{bmatrix}
			y[0] \\
			y[1] \\
			y[2]  \\
			\vdots \\
			y[n-1]  \end{bmatrix} =
		\begin{bmatrix}
			h[0] & h[n-1] &h[n_t-2] & \dots & h[1] \\
			h[1] & h[0] &h[n_t-1] & \dots & h[2] \\
			h[2] & h[1] &h[0] & \ddots & \vdots \\
			\vdots & h[2] &h[1] & \ddots & h[n-1] \\
			h[n-1] & \dots &h[2] & h[1] & h[0]   \end{bmatrix}  &\begin{bmatrix}
			u[0] \\
			u[1] \\
			u[2]  \\
			\vdots \\
			u[n-1]  \end{bmatrix} \,.
	\end{split}
\end{equation}

On the other hand, we know from the convolution theorem that the same can be achieved in the frequency domain using an entry by entry multiplication of the Fourier transform of the signals and the one of the impulse response. In terms of matrix multiplications, these are, respectively:

\begin{equation}
	\label{s8_eq24}
	\widehat{\mathbf{u}}=\overline{\bm{\Psi}}_{\mathcal{F}} \mathbf{u}\,, \quad \widehat{\mathbf{y}}=\overline{\bm{\Psi}}_{\mathcal{F}} \mathbf{y}\,,\mbox{and}\,\, 
	\widehat{\mathbf{H}}=\overline{\bm{\Psi}}_{\mathcal{F}} \mathbf{h}\,
\end{equation} where $\mathbf{h}:= h[k]\in \mathbb{R}^{n_t\times 1}$ is the impulse response arranged as a column vector. To perform the entry-by entry matrix multiplication, we define a diagonal matrix $\bm{H}=diag(\widehat{\mathbf{H}})$. Hence the operation in the frequency domain becomes $\widehat{\mathbf{y}}=\bm{H}\,\widehat{\mathbf{u}}$. Introducing the Fourier transforms in \eqref{s8_eq24} and moving back to the time domain we get

\begin{equation}
	\label{s8_eq25}
	\widehat{\mathbf{y}}=\bm{H}\,\widehat{\mathbf{u}}\rightarrow \overline{\bm{\Psi}}_{\mathcal{F}} \mathbf{y}=\bm{H}\,\overline{\bm{\Psi}}_{\mathcal{F}} \mathbf{u}\rightarrow {\mathbf{y}}=\bm{\Psi}_{\mathcal{F}} \bm{H}\,\overline{\bm{\Psi}}_{\mathcal{F}}\mathbf{u}\quad \mbox{hence} \,\, \bm{C}_h=\bm{\Psi}_{\mathcal{F}} \bm{H}\,\overline{\bm{\Psi}}_{\mathcal{F}}
\end{equation} The last equation on the right is obtained by direct comparison with \eqref{s8_eq24} and is extremely important: this is an eigenvalue decomposition. Each $\widehat{\mathbf{H}}[n]$ entry of the frequency response vector is an eigenvalue of $\bm{C}_h$ and the Fourier basis element $\bm{\psi}_{\mathcal{F}n}$ (corresponding to the frequency $f_n$) is the associated eigenvector.

Of great interest is the case in which the impulse response is even, i.e. $h[k]=h[-k]$. In this case the circulant matrix $\bm{C}_h$ becomes also symmetric and its eigenvalues are real: this is the case of zero-phase (non-causal) filters. In this special case, also the eigenvectors can be taken as real harmonics: they could be either sinusoidals or cosinusoidals, i.e. the basis of the Discrete Sine Transform (DST) or the Discrete Cosine Transform (DCT) \citep{Strang2007}.
In a completely different context, this eigenvalue decomposition also appear in the Proper Orthogonal Decomposition of a special class of signals. Because of its importance, the reader is encouraged to pause and practice with the following exercise.

\begin{tcolorbox}[breakable, opacityframe=.1, title=Exercise 1: The DFT and the Diagonalization of Circulant Matrices]
	
	Consider the discrete signal $\mathbf{u}$ from Exercise 5 in Chapter 4. 1) compute the DFT via Matrix multiplication and compare that the result with the one from numpy's \textsc{fft} routine. 2) Prepare the impulse response $h[k]$ of a low pass filter of order $N=211$ using the windowing method and a Hamming window. Then implement this filter using the matrix multiplication form of the convolution in \eqref{s8_eq23} and the matrix multiplication in the frequency domain in \eqref{s8_eq25}. Compare the results with any of the previous implementation from chapter 4. 3) Construct the convolution matrix from \eqref{s8_eq25} and test the validity of the eigenvalue decomposition previously introduced. 
	
\end{tcolorbox}

\subsection{The Proper Orthogonal Decomposition (POD)}\label{sec_8_3_3}

The Proper Orthogonal Decomposition (POD) is introduced in Chapter 6 we here only focus on its special properties in relation to the algorithm 2.

The temporal structure of the POD are eigenvectors of the temporal correlation matrix, which in the notation of this chapter reads

\begin{equation}
	\label{s8_eq26}
	\bm{K}=\bm{D}^{\dag} \bm{D}= \bm{\Psi}_{\mathcal{P}}  \bm{\Lambda} \bm{\Psi}^T_{\mathcal{P}}\quad \mbox{i.e.}\quad \bm{K}[i,j]=\sum^{n_t-1}_{r=0}\sigma^2_{\mathcal{P}}\,\bm{\psi}_{\mathcal{P}}[i] \bm{\psi}_{\mathcal{P}}[j]
\end{equation}

The subscript $\mathcal{P}$ is used to distinguish the POD.
The notation on the right is based on an outer product representation of the eigenvalue decomposition of symmetric matrices; this will be useful in Section \ref{sec_8_4}.
The first key feature of the decomposition is that the eigenvalues of $\bm{K}$ are linked to the POD amplitudes as $\bm{\Lambda}=\bm{\Sigma}^2_{\mathcal{P}}$. Hence the diagonalization in \eqref{s8_eq26} \emph{also provides the POD amplitudes} and the normalization in line 4 of the algorithm 2 is not needed. Introducing $\bm{\Psi}_{\mathcal{P}}$ in this algorithm gives the Sirovinch's formulation\footnote{actually a much less efficient version: the projection step in line 1 of the algorithm could simply be $\tilde{\bm{D}}=\bm{D}\bm{\Psi}_{\mathcal{P}}$} described in Chapter 6. Observe that since the POD amplitudes are the singular values of the dataset matrix, it is possible to compute the convergence in \eqref{s8_eq19} without computing norms:

\begin{equation}
	\label{s8_eq27}
	E(\tilde{r})=\frac{||\bm{D}-\tilde{\bm{D}}(\tilde{r})||_2}{||\bm{D}||_2}=\sqrt{\frac{\sum^{R-1}_{r=\tilde{r}}\sigma^2_{\mathcal{P}r}}{\sum^{R-1}_{r=0}\sigma^2_{\mathcal{P}r}}}
\end{equation}

For the following discussion, two observations are of interest.
The first is that introducing the POD factorization (i.e., the SVD) in \eqref{s8_eq26} we have

\begin{equation}
	\label{s8_eq28}
	\bm{K}=\bm{\Psi}_{\mathcal{P}}\bm{\Sigma}^{-1}_{\mathcal{P}}\bm{\Phi}^T_{\mathcal{P}}\bm{\Phi}_{\mathcal{P}}\bm{\Sigma}^{-1}_{\mathcal{P}}\bm{\Psi}^T_{\mathcal{P}}\quad \mbox{i.e}\quad \bm{\Lambda}=\bm{\Sigma}^{-1}_{\mathcal{P}}\bm{\Phi}^T_{\mathcal{P}}\bm{\Phi}_{\mathcal{P}}\bm{\Sigma}^{-1}_{\mathcal{P}}\,.
\end{equation}

The last step arise from a direct comparison with \eqref{s8_eq26}. Because $\bm{\Lambda}$ is diagonal, we see that we must have $\bm{\Phi}^T_{\mathcal{P}}\bm{\Phi}_{\mathcal{P}}=\mathbf{I}$, i.e. the spatial structures \emph{are also orthonormal}. We know this from Chapter 6, where it was shown that these are eigenvectors of the spatial correlation matrix. The key observation is that the reverse must also be true: \emph{every decomposition that has orthonormal temporal and spatial structure is a POD}. This brings to the second observation on the uniqueness of the POD. In a dataset that leads to modes of equal energetic importance, the amplitudes of the associated POD modes tend to be equal. This means repeated eigenvalues of $\bm{K}$ and thus non-unique POD. In the extreme case of a purely random dataset, it is easy to see that the POD modes are all equal \citep{Mendez2017}, and there are infinite possible PODs. The impact of noise in a POD decomposition is further discussed in Chapter 9.

Finally, observe that the POD is based on error minimization (or, equivalently, amplitude maximization) and has no constraints on the frequency content in its temporal structures $\bm{\Psi}_{\mathcal{P}}$. However, a special case occurs in an ideally stationary process. In such a process, the temporal correlations are invariant with respect to time delays and solely depend on the time lag considered in the correlation. Hence the correlation $\bm{K}[1,4]=\mathbf{d}^T_1\,\mathbf{d}_4$ is equal to  $\bm{K}[4,7]=\mathbf{d}^T_4\,\mathbf{d}_7$ or $\bm{K}[11,14]=\mathbf{d}^T_{11}\,\mathbf{d}_{14}$, for example. 

In other words, the temporal correlation matrix $\bm{K}$ becomes circulant, like the matrix $\mathbf{C}_h$ in \eqref{s8_eq23} for the convolution. Therefore, its eigenvectors are harmonics: we conclude that \emph{the POD of an ideally stationary dataset is a either a DCT or a DST\footnote{depending on whether the temporal average has been removed: the mean is accounted for in a DCT but not in a DST.}}. This property is the essence of the Spectral POD proposed by \cite{sieber2016spectral}, which introduces an ingenious FIR filter along the diagonal of $\bm{K}$ to reach a compromise between the energy optimality of the POD modes and the spectral purity of Fourier modes. Depending on the strength of this filter, the SPOD offers an important bridge between the two decompositions.

\subsection{The Dynamic Mode Decomposition (DMD)}\label{sec_8_3_4}

The Dynamic Mode Decomposition (DMD) is introduced in Chapter 7 and we here focus on its link with the algorithm 2. The DMD aims at fitting a linear dynamical system to the data, writing every snapshot as 

\begin{equation}
	\label{s8_eq29}
	\mathbf{d}_{\mathbf{i}}[k+1]=\sum^{n_s-1}_{r=0} a_r \bm{\psi}_{\mathcal{D}r}[\mathbf{i}] e^{-p_r t_{k}}=\sum^{n_s-1}_{r=0} a_r \bm{\psi}_{\mathcal{D}r}[\mathbf{i}] \lambda_r^{k-1}\,.
\end{equation}

The subscript $\mathcal{D}$ is used to distinguish the DMD.
In Chapter 4 we have seen that a complex exponentials are eigenfunctions of linear dynamical systems; Chapter 10 describes their state-space representation. Observe that the DMD is a sort of Z-transform in which the basis elements only include the \emph{poles} of the systems (see Chapter 4), i.e. the eigenvalues of the propagating matrix $\bm{A}$ in the state-space representation in Chapter 10. Such decomposition is natural for a homogeneous linear system, while the extension of the DMD for the forced system is proposed by \cite{Proctor2016}. This extension makes the DMD an extremely powerful system identification tool that can be combined with the linear control methods introduced in Chapter 10.

Many algorithms have been developed to compute the eigenvalues $\lambda_r$'s from the dataset (see Chapter 7). Once these are computed, the matrix containing the temporal structures of the DMD, can be constructed:

\begin{equation}
	\label{s8_eq30}
	\bm{\Psi}_{\mathcal{D}}=\begin{bmatrix}
		1& 1  & 1  & \dots &1   \\
		\lambda_1& \lambda_2 & \lambda_3 & \dots & \lambda_{n_t} \\
		\lambda^2_1 & \lambda^2_2 & \lambda^2_3  & \dots & \lambda^2_{n_t} \\
		\vdots & \vdots &  \ddots  & \vdots \\
		\lambda^{n_t}_1&\lambda^{n_t}_2 & \lambda^{n_t}_3& \dots & \lambda^{n_t}_{n_t} \\
	\end{bmatrix}\in \mathbb{C}^{n_t\times n_t}\,\,,
\end{equation} and the decomposition completed following algorithm 2.

From a dynamical system perspective, we know from Chapter 4 that a system is stable if all its poles are within the unit circle, hence if all the $\lambda$'s have modulus $|\lambda|\leq 1$. From an algebraic point of view, we note that $\bm{\Psi}_{\mathcal{D}}$ differs from the temporal structures of the other decompositions in two important ways. First, it is generally not orthonormal: the full projection must be considered in line 1 of algorithm 2. Second, its inverse might not exist: convergence is not guaranteed. 

Different communities have different ways of dealing with this lack of convergence. In the fluid dynamics community, more advanced DMD algorithms such as the sparsity promoting DMD \cite{Jovanovic2014pof} or the optimized DMD \cite{Chen2012}, have been developed to enforce that all the $\lambda$'s have unitary modulus. This is done by introducing an optimization problem in which the cost function is defined on the error minimization of the DMD approximation. Consequently, vanishing or diverging decompositions are penalized. Observe that having all the modes \emph{on the unit circle does not necessarily imply that these are DFT modes}: the DMD does not impose any orthogonality condition, which would force the modes to have frequencies that are multiples of a fundamental one. If orthogonality is enforced as a constraint to the optimization, then the only degree of freedom distinguishing DMD and DFT is in the choice of the fundamental tone: while this is $T=n_t\Delta t$ for the DFT, the DMD can choose different values and by-pass problems like spectral leakage or windowing \cite{WINDOWING}.

In the climatology community, in which the DMD is known as Principal Oscillation Patterns  \citep{Hasselmann1988,Storch1990}{POP, see} or Linear Inverse Modeling  \citep{Penland1996,LIM1}{LIM, see}, the lack of convergence is seen as a natural limit of the linearization process. Accordingly, the results of these decompositions are usually presented in terms of \emph{e-folding time}, which is the time interval within which the temporal structures decrease or increase by a factor $e$. In time-series analysis, modes with a short e-folding time are modes with shorter predictive capabilities. Variants of the POP/LIM have been presented in the fluid dynamics community under the name of Oscillating Pattern Decomposition  \citep{Vlacav1}{OPD, see}. The link between DMD, POP, and LIM is also discussed by \cite{Schmid2010jfm} and \cite{Tu2014jcd}.

\bigskip

\begin{tcolorbox}[breakable, opacityframe=.1, title=Exercise 2: Eigenfunction Expansion and Modal Analysis]
	
	Consider a 2D pulsating Poiseuille flow within two infinite plates at a distance $2b$ over a length $L$, as shown in Figure 2a. Assume that the flow is fully developed ($\partial_x u=0$) and incompressible ($\partial_x u+ \partial_y v=0$). It is easy to show that the pressure gradient must thus be uniform ($\partial_x p=\mbox{const}$). If a sinusoidal pressure pulsation is introduced at the inlet, the Navier Stokes equation governing the velocity profile in the stream-wise direction are
	
	\begin{equation}
		\label{s8_eq31}
		\rho 	\partial_t u=- \frac{\Delta p (t)}{L}+\mu \partial_{yy} u \quad \mbox{with} \quad \Delta p(t)=p_M+p_A \cos(\omega t)\,
	\end{equation} where $\rho$ and $\mu$ are the density and the dynamic viscosity of the fluid. Introducing the following reference quantities
	
	\begin{equation}
		\label{s8_eq32}
		[y]=b \quad [x]=L \quad [p]=p_M \quad [t]=1/\omega \,,
	\end{equation} the problem can be written in dimensionless form as follows:

	\begin{equation}
		\label{s8_eq33}	
		\begin{dcases} 
			\partial_t \hat{u} -\frac{1}{\mathcal{W}^2} \partial_{\hat y\hat y}  \hat{u}=-\frac{\hat{p}_A}{\mathcal{W}^2}  \cos(\hat t) & -1<\hat{y}<1, \, t\geq0\\
			\hat{u}(-1,\hat{t})=\hat{u}(1,\hat{t}) = 0 & t\geq0\\
			\hat{u}(\hat{y},0)=u_0 (\hat{y}) & -1<\hat{y}<1 \,\,.
		\end{dcases}
	\end{equation} Here the hat denotes dimensionless variables, i.e $\hat{a}=a/[a]$ is the dimensionless scaling of a variable $a$ with respect to the reference quantity $[a]$. The dimensionless number controlling the response of the profile is the \emph{Womersley} number $\mathcal{W}=b\sqrt{\omega \rho/\mu}$.
	
	Eq \eqref{s8_eq33} is a linear parabolic PDE with homogeneous Dirichlet boundary conditions on both sides of the domain. A variable separated solution can be found by expanding the velocity profile in terms of eigenfunctions of the Laplacian operator. Given the boundary conditions in \eqref{s8_eq32}, these are cosine functions. The solution via eigenfunction expansion \citep{Mendez_POI}:
	
	\begin{equation}
		\label{EIG_SOLUTION}
		\hat u(\hat y, \hat t)=\sum_{n=1}^{\infty}\phi_{\mathcal{E}n}(\hat{y})\, \sigma_{\mathcal{E}n} \,  \psi_{\mathcal{E}n}(\hat t) 	
	\end{equation}

	with
	
	\begin{equation}
		\label{EIG_SOLUTION2}
		\begin{split}
			\phi_{\mathcal{E}n}&=\cos\biggl[(2n-1)\frac \pi 2 \hat y \bigg]\\
			\sigma_{\mathcal{E}n}&=\frac{16 \hat{p}_A}{(2n-1)\pi \sqrt{16 \mathcal{W}^4+(2n-1)^4\pi^4}}\\
			\psi_{\mathcal{E}n}&=(-1)^n\cos \biggl [\hat{t}-\tan^{-1}\biggl(\frac {\mathcal{W}}{[(2n-1)\pi]}\biggr)^2\biggr ]
		\end{split}
	\end{equation}
	
	Note that the spatial structures are orthogonal but the temporal ones are not. The amplitudes of the modes decay as $\propto1/(2n-1)^3$ when $\mathcal{W}\rightarrow 0$ and as $\propto1/(2n-1)$ when $\mathcal{W}\rightarrow \infty$. Finally, note that all the amplitudes tend to $\sigma_{\mathcal{E}n}\rightarrow 0$ as $\mathcal{W}\rightarrow \infty$: as the frequency of the perturbation increases, the oscillations in the velocity profile are attenuated. Consider a case with $\mathcal{W}=10$ and $\hat{p}_a=60$. Several dimensionless velocity profiles during the oscillations are shown in Figure 2b.

	In this exercise, the reader should compare the DFT, the DMD, and the POD with the eigenfunction solution. Assume that the space discretization consists of $n_y=2000$ while the time discretization consists of $n_t=200$ with a dimensionless sampling frequency of $\hat{f}_s=10$. 
	
	First, construct the discrete dataset from \eqref{EIG_SOLUTION}-\eqref{EIG_SOLUTION2} by setting it in terms of the canonical factorization in \eqref{s8_eq11}. Then, prepare a function that implements the algorithm 2 described in the previous section and use this algorithm to compute the DFT, POD, and DMD. Finally, plot the amplitude decay of all the decompositions and show the first three dominant structures in space and time for each.
	
	\begin{center}%

		\begin{minipage}{.49\linewidth}
			\centering
			\includegraphics[width=0.8\linewidth]{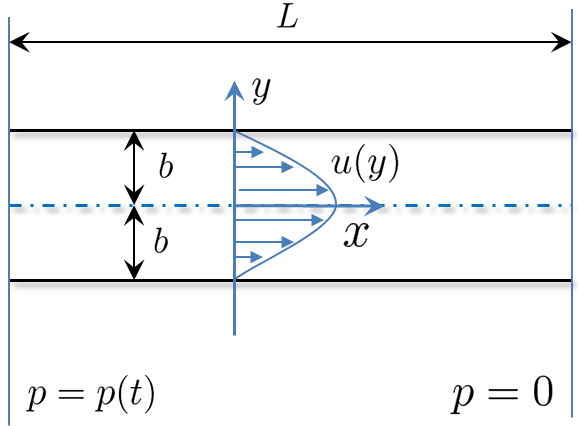}%
			\label{a)}%
		\end{minipage}
		\begin{minipage}{.49\linewidth}
			\includegraphics[width=\linewidth]{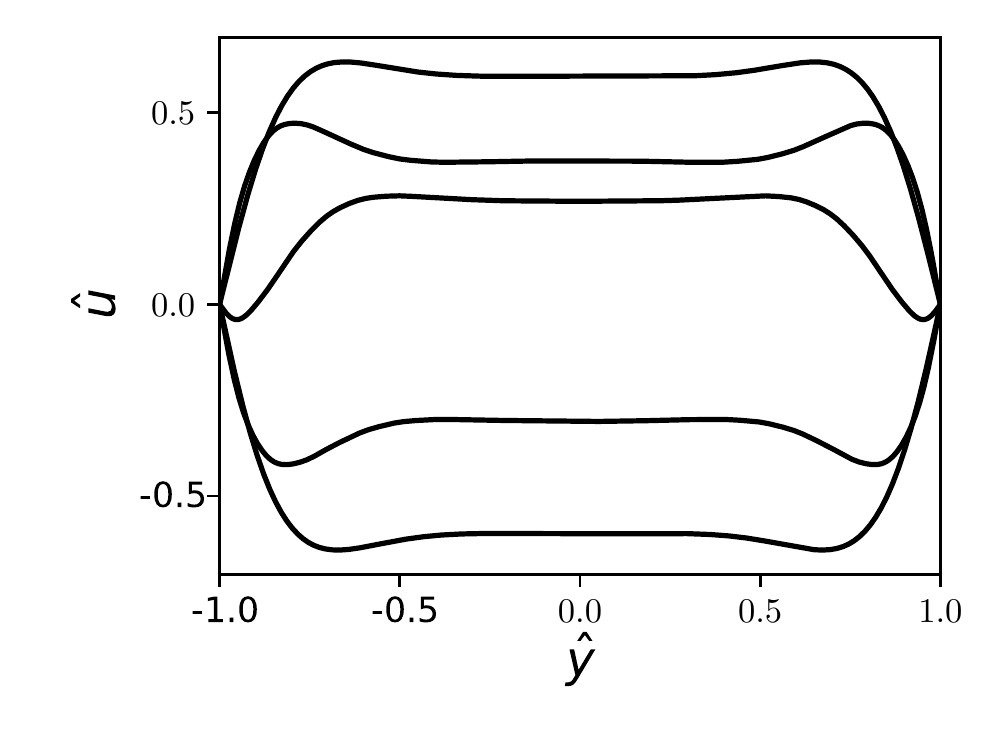}%
			\label{b)}%
		\end{minipage}			
		{Fig (a): Schematic of the 2D pulsating Poiseuille Flow considered in this exercise. Fig (b) snapshots of the velocity profile during the pulsation.}	
		\label{ch8fig2}
	\end{center}

\end{tcolorbox}

\section{The Multiscale POD (mPOD) }\label{sec_8_4}

This section is composed of three parts. We first analyze in Section \ref{sec_8_4_1} the impact of frequency constraints on the POD. In particular, we are interested in what happens if we filter a dataset before computing the POD. We will see that if a frequency is removed from the data, it is also removed from the temporal structures of its POD. We extend this result to a decomposition of the data via MRA (see Chapters 4 and 5). The MRA uses a filter bank to break the data into \emph{scales} and we here see how to keep the PODs of all scales \emph{mutually orthogonal}. Then, the POD bases of each scale can be assembled into a single orthonormal basis. That is the basis of the Multiscale Proper Orthogonal Decomposition-- the mPOD, presented in Section \ref{sec_8_4_2}. Finally, the mPOD algorithm is described in section \ref{sec_8_4_3}.

\subsection{Frequency Constrained POD}\label{sec_8_4_1}

We compute the POD of a dataset that has first been filtered. \emph{Are the POD modes also filtered?} \emph{can we obtain the same modes by filtering the data or filtering the temporal correlation matrix?} \emph{how are these linked to the POD modes of the original data?} This section answers these questions. 

First, we need some definitions. Let the filter have a finite impulse response\footnote{Which thus implies that the filter cannot be `ideal', as discussed in Chapter 4.} $h\in \mathbb{R}^{n_t\times 1}$. Following Section \ref{sec_8_3_2}, the frequency response in the discrete domain is a complex vector $\widehat{\mathbf{H}}=\bm{\Psi}_{\mathcal{F}} \bm{h}\in\mathbb{C}^{n_t\times 1}$. This vector can be arranged in two possible matrices, both useful for the following discussion. The first is linked to the diagonalization of circulant matrices, and is introduced in \eqref{s8_eq25}. The second consists of a row-wise replication of the same vector. These are $\bm{H}\in\mathbb{C}^{n_t\times n_t}$ and $\bm{\mathcal{H}}=\mathbf{1}\cdot\widehat{\mathbf{H}}^T\in\mathbb{C}^{n_s\times n_t}$, with $\mathbf{1}\in\mathbb{R}^{n_s\times1}$ a column vector of ones\footnote{To be normalized by $1/\sqrt{n_s}$ to keep it of unitary length.}. To avoid confusion, we here show both matrices\footnote{Note that the right multiplication by $\bm{H}$ is equivalent to the Hadamard product by $\bm{\mathcal{H}}$} :

\begin{equation}
	\label{s8_eq36}
	\bm{H}=\begin{bmatrix}
		\widehat{\mathbf{H}}[0]& 0 & 0 &  0\\
		0& \widehat{\mathbf{H}}[1] & 0 & 0\\
		\vdots & \vdots &  \ddots  & \vdots \\
		0& 0 & 0 & \widehat{\mathbf{H}}[n_t] \\
	\end{bmatrix}\quad \quad \bm{\mathcal{H}}=\begin{bmatrix}
		\widehat{\mathbf{H}}[0]& \widehat{\mathbf{H}}[1] & \dots &  \widehat{\mathbf{H}}[n_t]\\
		\widehat{\mathbf{H}}[0]& \widehat{\mathbf{H}}[1] & \dots &  \widehat{\mathbf{H}}[n_t]\\
		\vdots & \vdots &  \vdots  & \vdots \\
		\widehat{\mathbf{H}}[0]& \widehat{\mathbf{H}}[1] & \dots &  \widehat{\mathbf{H}}[n_t] \\
	\end{bmatrix}
\end{equation}

Recall that $\bm{H}$ is the eigenvalue matrix of the circulant convolution matrix $\bm{C}_h$ in \eqref{s8_eq23}. Assuming that the filtering in the time domain is performed in all the spatial location equally, the associated matrix multiplications are:

\begin{equation}
	\label{s8_eq37}
	\bm{D}_{\mathcal{H}}=\bm{D}\,\bm{C}_h=\overbrace{\underbrace{\bm{D} \overline{\bm{\Psi}}_{\mathcal{F}}}_{\widehat{\bm{D}}}\bm{H} }^{\widehat{\bm{D}}_{\mathcal{H}}}\bm{\Psi}_{\mathcal{F}} \quad \mbox{or}\quad 
	\bm{D}_{\mathcal{H}}=\overbrace{\bigr[\underbrace{\bigl(\bm{D}\, \overline{\bm{\Psi}}_\mathcal{F}\bigr)}_{\widehat{\bm{D}}}  \odot\bm{\mathcal{H}}\bigr]}^{\widehat{\bm{D}}_{\mathcal{H}}}{\bm{\Psi}}_\mathcal{F}\,\,
\end{equation} where $\odot$ is the Hadamard product, i.e. the entry by entry multiplication between two matrices. Observe that $\widehat{\bm{D}}_{\mathcal{H}} \bm{\Psi}_{\mathcal{F}}$ is the inverse Fourier transform of $\widehat{\bm{D}}_{\mathcal{H}}$ (i.e. the frequency spectra of the filtered data) while $\widehat{\bm{D}}=\bm{D}\bm{\Psi}_{\mathcal{F}}$ is the Fourier transform of $\bm{D}$ along its rows (i.e. in the time domain).

The representation on the left is a direct consequence of the link between convolution theorem and eigenvalue decomposition of circulant matrices introduced in \eqref{s8_eq25}. The representation on the right opens to an intuitive and graphical representation of the filtering process that is worth discussing briefly.  

First, we introduce the cross-spectral density matrix $\bm{K}_{\mathcal{F}}$. This collects the inner product between the frequency spectra of the data evolution; it is the frequency-domain analogous of the temporal correlation matrix $\bm{K}$. These matrices are linked: 

\begin{equation}
	\label{s8_eq38}
	\bm{K}_{\mathcal{F}}=\widehat{\bm{D}}^\dag\,\widehat{\bm{D}}= \bm{\Psi}_{\mathcal{F}} \,  \bigl [\bm{D}^\dag\, \bm{D}\bigr ] \,\overline{\bm{\Psi}}_{\mathcal{F}}=\bm{\Psi}_{\mathcal{F}} \,  \bm{K} \,\overline{\bm{\Psi}}_{\mathcal{F}} \Longleftrightarrow  \bm{K}=\overline{\bm{\Psi}}_{\mathcal{F}} \,  \bm{K}_{\mathcal{F}} \,{\bm{\Psi}}_{\mathcal{F}}\,.
\end{equation}

Since $\bm{\Psi}_{\mathcal{F}}\overline{\bm{\Psi}}_{\mathcal{F}}=\mathbf{I}$, the equation on the right is a similarity transform, hence $\bm{K}_{\mathcal{F}}$ and $\bm{K}$ share the same eigenvalues.  Introduce the eigenvalue decomposition of $\bm{K}$ in \eqref{s8_eq26}:

\begin{equation}
	\label{s8_eq39}
	\bm{K}_{\mathcal{F}}=\bm{\Psi}_{\mathcal{F}} \,\bigl( \bm{\Psi}_{\mathcal{P}}  \bm{\Sigma}^2_{\mathcal{P}} \bm{\Psi}^T_{\mathcal{P}} \,\bigr)\overline{\bm{\Psi}}_{\mathcal{F}}= \overline{\widehat{\bm{\Psi}}}_{\mathcal{P}}\bm{\Sigma}^2_{\mathcal{P}} \overline{\widehat{\bm{\Psi}}}^{\dag}_{\mathcal{P}} \quad \mbox{i.e.}\quad \bm{K}_{\mathcal{F}}[i,j]=\sum^{n_t-1}_{r=0}\sigma^2_{\mathcal{P}}\,\overline{\widehat{\bm{\psi}}}_{\mathcal{P}}[i] {\widehat{\bm{\psi}}}_{\mathcal{P}}[j]
\end{equation}

We thus see that the eigenvectors of $\bm{K}_{\mathcal{F}}$ are the conjugate of the Fourier transform of the eigenvectors of $\bm{K}$. The outer product notation on the right shows that the diagonal of this matrix contains the sum of the power spectra of the temporal structures of all the modes. This is the sum of positive real quantities. 

Consider now the cross-spectral density matrix of the filtered data in \eqref{s8_eq36} and use the distributive property of the Hadamard product to get:

\begin{equation}
	\label{s8_eq40}
	\bm{K}_{\mathcal{F}}=\widehat{\bm{D}}_{\mathcal{H}}^\dag\,\widehat{\bm{D}}_{\mathcal{H}}=\bigl(\widehat{\bm{D}} \odot \bm{\mathcal{H}}\bigr)^\dag\bigl(\widehat{\bm{D}} \odot\bm{\mathcal{H}}\bigr)
	=\bigl(\widehat{\bm{D}}^\dag  \widehat{\bm{D}} \bigr) \odot\bigl(\bm{\mathcal{H}}^\dag \odot\bm{\mathcal{H}}\bigr)
	=\bm{K}_{\mathcal{F}} \odot\underline{\bm{\mathcal{H}}},
\end{equation} having introduced the 2D symmetric frequency transfer function $\underline{\bm{\mathcal{H}}}=\bm{\mathcal{H}}^{\dag} {\bm{\mathcal{H}}}$. Figure \ref{ch8fig2} gives a pictorial representation of the magnitude of such 2D frequency transfer function for ideal low-pass, band-pass and high-pass filters. In a low-pass filter, for example, the 2D frequency transfer function is unitary within a square $\bm{\mathcal{H}}[i,j]\neq 0 \,\forall i,j \in[-i_c,i_c]\times [-i_c,i_c]$ and approximately zero everywhere else.

We thus see that the filtering of the data constraints the spectral content of the POD modes: writing the cross-spectral density of the filtered data in terms of its eigenvectors as in \eqref{s8_eq38}, we see: 

\begin{equation}
	\label{s8_eq41}
	\bm{K}_{\mathcal{F}\mathcal{H}}=\bm{K}_{\mathcal{F}} \odot\underline{\bm{\mathcal{H}}}\Longleftrightarrow \bm{K}_{\mathcal{F}\mathcal{H}}[i,j]=\sum^{n_t-1}_{r=0}\sigma^2_{\mathcal{P}\mathcal{H}}\,\overline{\widehat{\bm{\psi}}}_{\mathcal{P}\mathcal{H}}[i] {\widehat{\bm{\psi}}}_{\mathcal{P}\mathcal{H}}[j]
\end{equation}

\begin{figure}[t]
	\centering
	Low Pass Filtering\\
	\includegraphics[keepaspectratio=true,width=0.9
	\columnwidth]{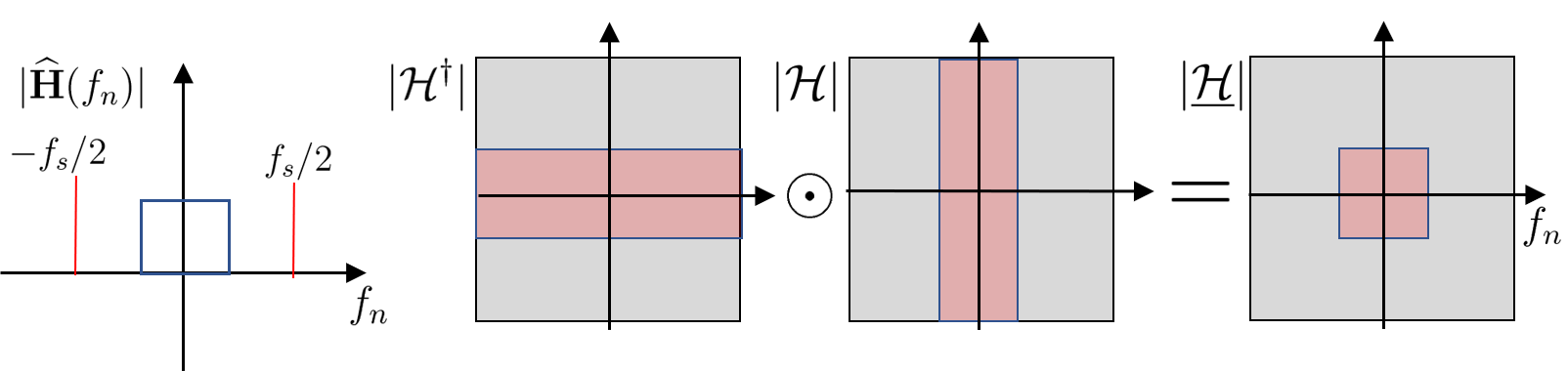}\\
	\vspace{2mm}
	Band Pass Filtering\\
	\includegraphics[keepaspectratio=true,width=0.9 \columnwidth]{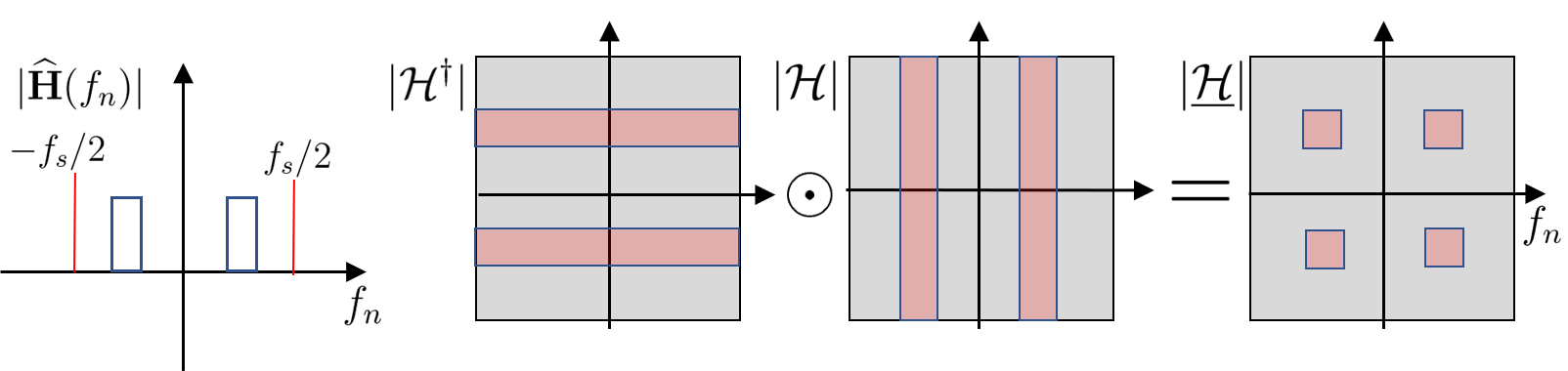}\\
	High Pass Filtering\\
	\includegraphics[keepaspectratio=true,width=0.9 \columnwidth]{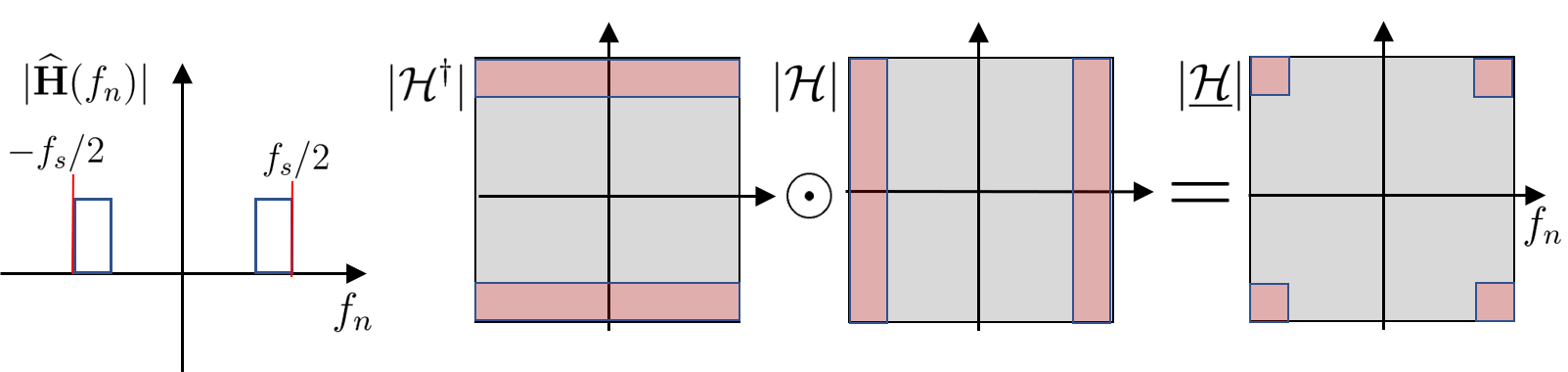}
	\caption{Pictorial representation of the magnitudes of the frequency transfer functions of the filters considered in this section. The figure on the left shows the magnitude of the 1D frequency response $\mathbf{\widehat{H}}$. On the right, the figure shows how 2D filters $\underline{\bm{\mathcal{H}}}=\bm{\mathcal{H}}^{\dag}\bm{\mathcal{H}}$ are constructed from the extended response $\bm{\mathcal{H}}$.
		The red region represents the pass-band portion, i.e. within which $|\bm{\mathcal{H}}|\approx 1$ while in the gray area it is $|\bm{\mathcal{H}}|\approx 0$. }
	\label{ch8fig8}
\end{figure}

Since $\bm{K}_{\mathcal{F}\mathcal{H}}\approx 0$ outside the band-pass range of the 2D filter, and since the entries along its diagonals are a summation of positive quantities, we conclude that \emph{frequencies removed from the data cannot be present in any of the POD modes}. Moreover, we can compute the POD modes of the filtered data from filtered the cross-spectral matrix $\bm{K}_{\mathcal{F}\mathcal{H}}$ or, more conveniently, from the temporal correlation matrix $\bm{K}_{\mathcal{H}}$. The link between these matrices can be written as 

\begin{equation}
	\label{s8_eq42}
	\bm{K}_{\mathcal{H}}=\bm{D}^\dag_\mathcal{H}\bm{D}_\mathcal{H}=\overline{\bm{\Psi}}_{\mathcal{F}} \bm{K}_{\mathcal{F}\mathcal{H}}\bm{\Psi}_{\mathcal{F}}=\overline{\bm{\Psi}}_{\mathcal{F}} \widehat{\bm{K}}_{\mathcal{H}}\overline{\bm{\Psi}}_{\mathcal{F}} \Longleftrightarrow \widehat{\bm{K}}_{\mathcal{H}}=\bm{P}_{\pi}\bm{K}_{\mathcal{F}\mathcal{H}}\,,
\end{equation} where the definition of 2D Fourier transform in \eqref{s8_eq8} is introduced and $\bm{P}_{\pi}=\bm{\Psi}_\mathcal{F} \bm{\Psi}_\mathcal{F}$ is the permutation matrix obtained by multiplying by the Fourier matrix twice\footnote{A funny question after reading this: what happens if the Fourier transform is performed four times?}:

\begin{equation}
	\label{Permu}
	\bm{P}_\pi=\bm{\Psi}_\mathcal{F} \bm{\Psi}_\mathcal{F}=\overline{\bm{\Psi}}_\mathcal{F}\overline{\bm{\Psi}}_\mathcal{F}=\begin{bmatrix}
		1& 0 & \dots & 0&  0 \\
		0& 0 &  &  0&  1 \\
		0& 0 &  &  1& \vdots \\
		\vdots& \vdots  & \scalebox{-1}[1]{$\ddots$} &  & \vdots \\
		0& 1 & 0 & \dots &  0 \\
	\end{bmatrix}\,.
\end{equation}

Finally, to answer the last question from this subsection, use the multiplication by the diagonal matrix $\bm{H}$ in \eqref{s8_eq39} instead of the Hadamard product:

\begin{equation}
	\label{s8_eq44}
	\bm{K}_{\mathcal{F}\mathcal{H}}=\bigl(\widehat{\bm{D}} \bm{H} \bigr)^\dag\bigl(\widehat{\bm{D}} \bm{H}\bigr)
	=\bm{H}^{\dag} \widehat{\bm{D}}^{\dag} \bm{D}\bm{H}=\bm{H}^{\dag} \bm{K}_{\mathcal{F}}\bm{H}\,.
\end{equation}

Since $\bm{H}^{\dag}\bm{H}\neq \mathbf{I}$, this is not a similarity transform. Hence the eigenvalues of these two matrices are different: there is no simple link between $\sigma_{\mathcal{P}}$ and the $\sigma_{\mathcal{P}\mathcal{H}}$, nor between the POD of the data and the POD of the filtered data. We do not attempt to build such a connection: it suffices to notice that the POD of the filtered data is the optimal bases for the portion of data whose spectra fall within the band-pass region of the filter. In the terminology of MRA (see Chapter 4), we call this a \emph{scale}.

\vspace{5mm}
\begin{tcolorbox}[breakable, opacityframe=.1, title=Exercise 3: The POD of Filtered Data]
	
	Consider the test case from the previous exercise but assume that the forcing is composed of two terms. These are sinusoids associated to a Womersley number of $\mathcal{W}=1$ and $\mathcal{W}=4$, modulated by two windows We refer to these terms as $F_1$ and $F_2$, respectively. As the problem is linear, we can sum the response of the velocity profile to each term. The response of the velocity profile at the centerline, for each of the terms, are shown in Figure \ref{ch8fig6}.  
	
	In this exercise, we want to analyze the response of the velocity profile to the perturbation at the highest frequency. 1) Perform a standard POD and analyze the resulting structures. 2) Filter the data using a high-pass filter that isolate the forcing $F_2$ and compute the POD again. Compute this decomposition by filtering the data and by filtering the correlation matrix. Compare the results. As an optional exercise, perform DMD and DFT using the codes from the previous exercise: is it possible to separate the contribution of each forcing term with these decompositions? The answer is no. Show it!

	\begin{center}%
		\includegraphics[width=0.7\linewidth]{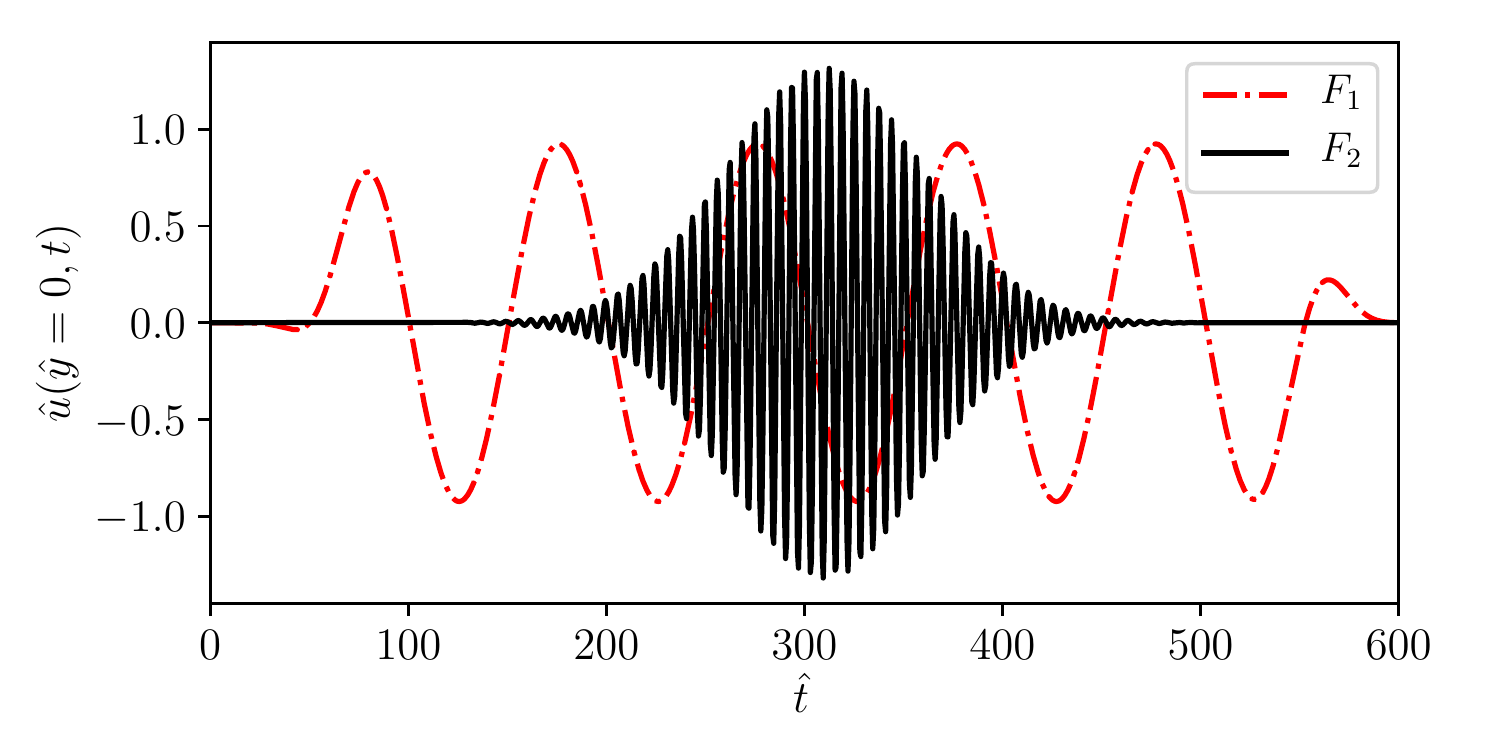}\\
		Time evolution of the velocity profile at the center of the channel when each of the two forcing terms $F_1$ and $F_2$ is active.
		\label{ch8fig6}
	\end{center}

\end{tcolorbox}

\medskip

\subsection{From Frequency Constraints to MRA }\label{sec_8_4_2}

The previous section provided the background to generalize the study of the frequency constrained POD to the Multiscale POD (mPOD). The mPOD constructs a set of orthogonal modes that have no frequency overlapping, that is no common frequency content within specific frequency bandwidths.

To illustrate the general idea, let us consider a dataset $\bm{D}\in \mathbb{R}^{n_s\times n_t}$, with $n_t<n_s$. We now use MRA to break the dataset into scales. With no loss of generalities, let us consider three scales; these are identified by three filters with  band-pass bandwidths $\Delta f_1=[0,f_1]$, $\Delta f_2=[f_1,f_2]$ and $\Delta f_3=[f_2, f_s/2]$. We lump this information in a frequency splitting vector $F_V=[f_1,f_2,f_3]$. With three scales, the temporal correlation matrix has nine partitions:

\begin{equation}
	\label{s8_eq45}
	\bm{K}=\bm{D}^{\dag} \bm{D}=\big(\bm{D}_1+\bm{D}_2+\bm{D}_3\bigr)^\dag \big(\bm{D}_1+\bm{D}_2+\bm{D}_3\bigr)=\sum^{3}_{i=1}\sum^{3}_{j=1} \bm{D}^{\dag}_i \bm{D}_j\,.
\end{equation}

Following a MRA formulation, we assume that these scales are isolated by filters with complementary frequency response, i.e.

\begin{equation}
	\label{s8_eq46}
	\bm{D}_{i}=\bigr[\bigl(\bm{D}\, \overline{\bm{\Psi}}_{\mathcal{F}}\bigr) \odot\bm{\mathcal{H}}_i\bigr]\bm{\Psi}_{\mathcal{F}}\quad \mbox{with} \quad\sum \bm{\mathcal{H}}_i\approx\mathbf{1} \quad \mbox{and} \quad \bm{\mathcal{H}}_i\odot\bm{\mathcal{H}}_j\approx\mathbf{0}\end{equation}

\begin{figure}[htbp]
	\centering
	\includegraphics[keepaspectratio=true,width=0.9 \columnwidth]{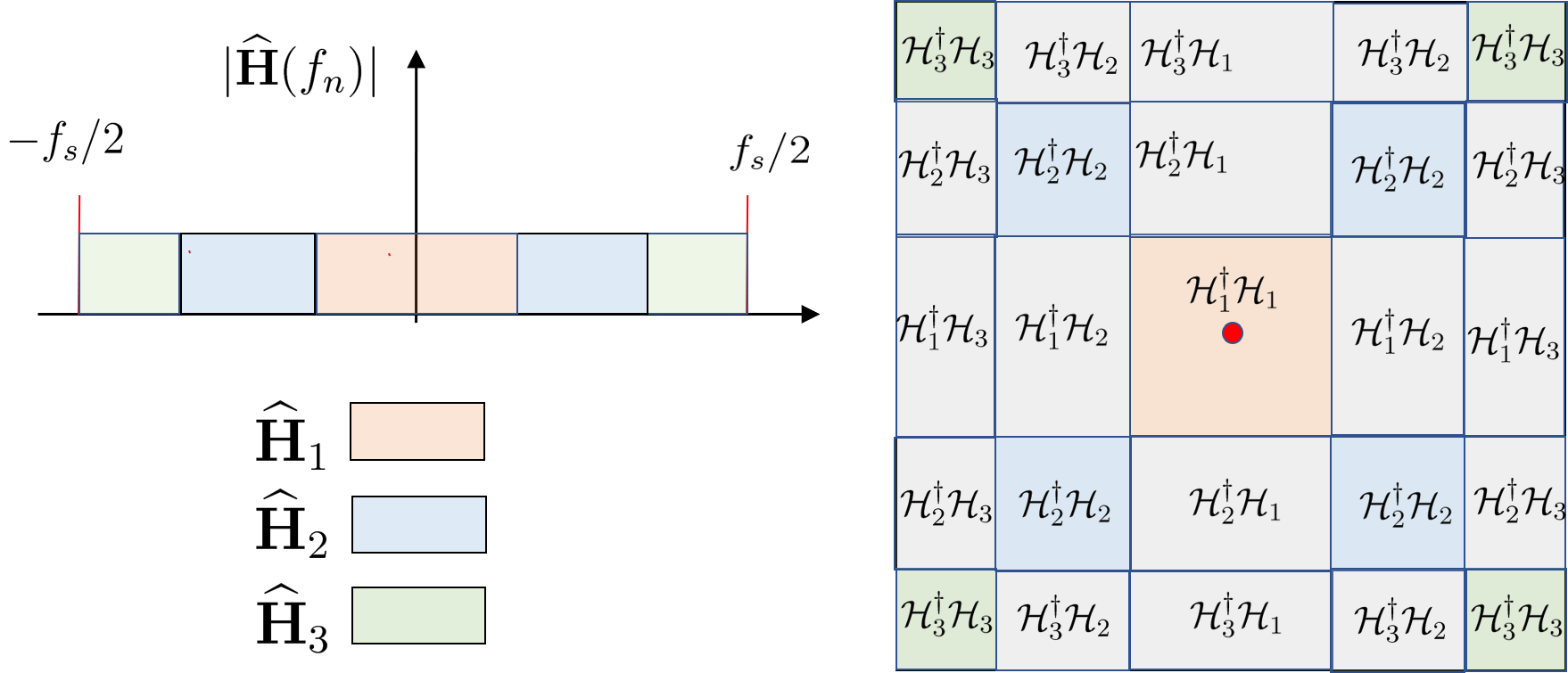}
	\caption[MRA Analysis and Cross-spectral Densities Repartition]{Repartition of the cross-spectral density matrix into the nine contributions identified by three scales. The origin is marked with a red circle. The `mixed terms' $\bm{\mathcal{H}}^{\dag}_i\bm{\mathcal{H}}_j$ with $i\neq j$ are colored in light gray, while the `pure terms' $\bm{\mathcal{H}}^{\dag}_i\bm{\mathcal{H}}_i$ are colored following the legend on the left, where the 1D transfer functions are also shown.}
	\label{ch8fig3}
\end{figure}

This ensures a lossless decomposition of the data. We now use \eqref{s8_eq39} and \eqref{s8_eq45} in \eqref{s8_eq44}, to analyze which portion of the cross-spectral density $\bm{K}_{\mathcal{F}}$ is taken by each of the contributions in \eqref{s8_eq44}. Figure \ref{ch8fig3} gives a pictorial view of the partitioning, which can be constructed following the graphical representation in Figure \ref{ch8fig2}. In Wavelet terminology, the term $\bm{\mathcal{H}}_1^\dag\bm{\mathcal{H}}_1$ is the \emph{approximation term} at the largest scale. The other `pure' terms $\bm{\mathcal{H}}_2^\dag\bm{\mathcal{H}}_2$ and $\bm{\mathcal{H}}_3^\dag\bm{\mathcal{H}}_3$ are \emph{diagonal details} of the scales $2$ and $3$. The terms $\bm{\mathcal{H}}_i^\dag\bm{\mathcal{H}}_j$ with $i>j$ are \emph{horizontal details} while those with $j<i$ are \emph{vertical details}.

We have seen in the previous section that removing a frequency from the dataset removes it from the temporal structures of its POD. There is thus no frequency overlapping\footnote{at least in case of ideal filters, which are perfectly complementary.} between the eigenvectors of the contributions $\bm{K}_1=\bm{D}_1^{\dag}\bm{D}_1$, $\bm{K}_2=\bm{D}_2^{\dag}\bm{D}_2$, $\bm{K}_3=\bm{D}_3^{\dag}\bm{D}_3$, which we denote as $\bm{K}_m=\bm{\Psi}_{\mathcal{H}m} \bm{\Lambda}_{\mathcal{H}m} \bm{\Psi}^{T}_{\mathcal{H}m}$. 
However, the same is not necessarily true for the correlations of the `mixed terms' $\bm{D}_i^{\dag}\bm{D}_j$ with $i\neq j$: the filters $\bm{\mathcal{H}}^{\dag}_i\bm{\mathcal{H}}_j$ with $i\neq j$ leave the diagonals of $\bm{K}_{\mathcal{F}}$ unaltered and thus the spectral constraints are much weaker: it is possible that eigenvectors from $\bm{D}_2^{\dag}\bm{D}_3$, for example, have frequencies that are already present among the eigenvectors of $\bm{D}_2^{\dag}\bm{D}_2$ or $\bm{D}^{\dag}_3\bm{D}_3$.

We opt for a drastic approach: remove all the `mixed terms' $\bm{D}^{\dag}_i\bm{D}_j$ with $i\neq j$ and consider the correlation as the sum of $M$ `pure terms'

\begin{equation}
	\label{s8_eq47}
	\bm{K}\approx \sum^{M}_{m=1} \bm{D}^{\dag}_{m}\bm{D}_m=\sum^{M}_{m=1} \bm{\Psi}_{\mathcal{F}}\bigr[\widehat{\bm{K}} \odot \bm{\mathcal{H}}^{\dag}_m \bm{\mathcal{H}_m}\bigr]\bm{\Psi}_{\mathcal{F}}=\sum^{M}_{m=1}\bm{\Psi}_{\mathcal{H}m} \bm{\Lambda}_{\mathcal{H}m} \bm{\Psi}^{T}_{\mathcal{H}m}\,,
\end{equation} having introduced \eqref{eq41} and the eigenvalue decomposition of each term.

This is a fundamental step to link the mPOD to POD and DFT. Such simplification results in a loss of information in the temporal correlation matrix, but it still ensures a complete reconstruction of the dataset. It remains to be seen if a basis constructed by combining all the eigenvectors of different scales can form a complete and orthonormal basis for $n_t$.

Consider first each of these scales independently. Let us assume that the number of non-zero frequencies in each scale, namely the number of frequency bins in the band-pass region, is $n_m$. If the frequency transfer functions are complementary, we have $\sum n_m=n_t$. If each scale contains at most $n_m$ frequencies, every temporal evolution of the data within that specific scale can be written as a linear combination of $n_m$ Fourier modes. Therefore, every contribution $\bm{D}_m$ is at most of $rank(\bm{D}_m)=n_m$ and hence has at most $n_m$ POD modes. These are are clearly orthogonormal, i.e. $\bm{\Psi}^{T}_{\mathcal{H}m}\bm{\Psi}_{\mathcal{H}m}=\mathbf{I}$ and form a basis for that specific scale.

Consider now modes from different scales, say for example $\bm{\psi}^{(1)}_{\mathcal{H}2}$ and $\bm{\psi}^{(4)}_{\mathcal{H}3}$, i.e. the first eigenvectors from scales $2$ and the fourth eigenvector from scale $3$. Because there is no frequency overlapping, their Fourier transforms are orthogonal, and thus they must also be orthogonal:

\begin{equation}
	\label{s8_eq48}
	\widehat{\bm{\psi}}^{(1)\dag}_{\mathcal{H}2}\widehat{\bm{\psi}}^{(4)}_{\mathcal{H}4}=0={\bm{\psi}}^{(1)\dag}_{\mathcal{H}2}\overline{\bm{\Psi}}_{\mathcal{F}}\bm{\Psi}_{\mathcal{F}}{\bm{\psi}}^{(4)}_{\mathcal{H}4}={\bm{\psi}}^{(1)\dag}_{\mathcal{H}2}{\bm{\psi}}^{(4)}_{\mathcal{H}4}
\end{equation}

We thus conclude that a complete and orthonormal basis for  $\mathbb{R}^{n_t}$ can be constructed from the eigenvectors of the various scales. That is the mPOD basis.

We close this section highlighting how the mPOD connects POD and DFT or DMD. In case of no frequency partitioning, the mPOD is a POD: the temporal structures can span the entire frequency range and are derived, as described in section \ref{sec_8_2_3} and Chapter 6, under the constraint of optimal approximation for any given number of modes. Introducing the frequency partitioning and the approximation in \eqref{s8_eq46}, we identify modes that are optimal only within the frequency bandwidth of each scale. 

\emph{The optimality of the full basis is lost as modes from different scales are not allowed to share the same frequency bandwidth}. As we introduce finer and finer partitioning, each mode is limited within a narrower frequency bandwidth. At the limit for $n_m=1$ and $M=n_t$, every mode is allowed to have only one frequency. The approximation in \eqref{s8_eq46} forces the spectra of the temporal correlation matrix to be diagonal, i.e., the correlation matrix is approximated as a Toeplitz circulant matrix. Accordingly, the mPOD tends towards the DFT or DMD depending on the boundary conditions used in the filtering of $\bm{K}$. If periodicity is assumed, the DFT is produced. If periodicity is not assumed, the decomposition selects harmonic modes whose frequency is not necessarily a multiple of the observation time, and a DMD is recovered.

\subsection{The mPOD Algorithm}\label{sec_8_4_3}

The pseudo-code for computing the the mPOD listed in algorithm 3. This is composed of eight steps.

\begin{algorithm}
	\KwIn{$\bm{D}\in \mathbb{R}^{n_s\times n_t}$, $F_V\in\mathbb{R}^{M\times 1}$}
	\KwOut{$\bm{\Psi}_{\mathcal{M}}\in\mathbb{R}^{n_t\times R}$, $\bm{\Sigma}_{\mathcal{M}}\in\mathbb{R}^{R\times R}$,  and $\bm{\Phi}_{\mathcal{M}}\in\mathbb{R}^{n_s\times R}$}
	Compute Temporal Correlation Matrix $\bm{K}=\bm{D}^{\dag}\bm{D}$
	
	Prepare the Filter bank $\bm{H}_1,\bm{H}_2 ... \bm{H}_M$

	Decompose $\bm{K}$ into $\bm{K}_1,\bm{K}_2... \bm{K}_m$ contributions 	
	
	Compute the temporal structures of each scale $\bm{K}=\bm{\Psi}_m\bm{\Lambda}_m\bm{\Psi}^T_m$
	
	Assembly basis
	$\bm{\Psi}^{0}=[\bm{\Psi}_1,... \bm{\Psi}_M]$ and $\bm{\lambda}^{0}=[diag(\bm{\Lambda}_1),... diag(\bm{\Lambda}_M)]$. 
	
	Prepare permutation $\bm{P}_{\Lambda}$, sorting $\bm{\lambda}^{0}$. Then shuffle columns $\bm{\Psi}^{1}= \bm{\Psi}^{0}\bm{P}_{\Lambda}$
	
	Enforce Orthogonality $\bm{\Psi}^{1}=\bm{\Psi}_{\mathcal{M}} \bm{R}\rightarrow \bm{\Psi}_{\mathcal{M}}=\bm{\Psi}^{1}\bm{R}^{-1}$.
	
	Complete the Decomposition using Algorithm 2
	$\bm{\Phi}_{\mathcal{M}}\rightarrow$ from $\bm{D}\,,\bm{\Psi}_{\mathcal{M}}$.
	\caption{Algorithm for computing the Multi-scale POD.}
	\label{algo:mPOD}
\end{algorithm}

The first step computes the temporal correlation matrix, as in the POD. The second step prepares the filter bank, according to the introduced frequency vector $F_V$. The third step computes the MRA of the temporal correlation matrix; the fourth step computes the eigenvectors of each\footnote{\emph{Can we by-pass this step? Can we compute all the $\bm{\Psi}_m$'s from one single diagonalization (of a properly filtered matrix?} These were two brilliant questions by Bo B. Watz, development engineer at Dantec Dynamics, who attended the lecture series. The answer is that it is possible, in principle, to compute the basis from one clever diagonalization. However, the filtering is challenging: that was the starting point of new exciting developments towards the \emph{fast mPOD}.}.

At this stage, we can proceed with the preparation of a single basis. First, in step $5$, all the temporal structures and associated eigenvalues are collected into a single matrix $\bm{\Psi}^{0}$ and a single vector of eigenvalues $\bm{\lambda}^{0}$. These eigenvalues give an estimation of the relative importance of each term. This is sorted in descending order and the information used to permute the columns of $\bm{\Psi}^{0}$ using an appropriate permutation matrix $\bm{P}_{\Lambda}$. 

If the filters were ideal, the resulting temporal basis $\bm{\Psi}^{1}= \bm{\Psi}^{0}\bm{P}_{\Lambda}$ would be completed. However, to compensate for the non-ideal frequency response of the filters, a reduced \textsc{QR} factorization is used in step seven to enforce the orthonormality of the mPOD basis. The result of the orthogonalization procedure is the matrix of the temporal structure of the mPOD. The final step is the projection to compute the spatial structures that is common to every decomposition and which can be computed using\footnote{Observe that a simplified version could be used since the temporal structure is orthonormal!} algorithm 2.

Note that special attention should be given to the filtering process in the third step. If this is done as in Exercise 4, it is implicitly assumed that the matrix $\bm{K}$ is periodic. When such an assumption is incorrect, edge effects could appear after the filtering process. The cure for these effects is described in classical textbooks on image processing \citep{Gonzalez2017} and include \emph{reflections} or \emph{extrapolations}.
To limit the scope of this chapter, we do not address these methods here\footnote{interested readers can find the related codes at \url{https://github.com/mendezVKI/MODULO}}. The reader should nevertheless be able to construct the mPOD algorithm using all the codes developed for the previous exercises. 

\medskip

\begin{tcolorbox}[breakable, opacityframe=.1, title=Exercise 4: Your own mPOD algorithm]
	
	Combine the codes from the previous exercises to build your own function for computing the mPOD, following algorithm 3. This function should compute $\bm{\Psi}_{\mathcal{M}}$, while the decomposition can be completed using algorithm 2. 
	
	Assume that filters are constructed using Hamming windows, and the filter order is an input parameter. Test your function with the previous exercise and show the structures of the first two mPOD modes. Compare their amplitudes with the ones of the POD modes.
	
\end{tcolorbox}

\section{Tutorial Test Cases}\label{sec_8_5}

In addition to the coding exercises, this chapter include two tutorial test cases from Time-Resolved Particle Image Velocimetry measurements. The first dataset collects the velocity field of a planar impinging gas jet in stationary conditions; the second is the flow past a cylinder in transient conditions. These are described in \cite{mendez2019multi} and \cite{DANTEC_VKI}. Other examples of applications of the mPOD on experimental data can be found in \cite{Mendez_Journal_2}, \cite{Mendez20192} and \cite{CLAU}.
More exercises and related codes can be found in \url{https://github.com/mendezVKI/MODULO}, together with an executable with Graphical User Interface (GUI) developed by \cite{MODULO}.

\subsection{Test Case 1}\label{tut1}

To download the dataset and store it into a local directory, run the following script

\begin{centering}
	\begin{lstlisting}[language=Python,linewidth=12cm,xleftmargin=0\textwidth,xrightmargin=0\textwidth,backgroundcolor=\color{yellow!10}]
import urllib.request
print('Downloading Data for Tutorial 1...')
url = 'https://osf.io/c28de/download'
urllib.request.urlretrieve(url, 'Ex_4_TR_PIV_Jet.zip')
print('Download Completed! I prepare data Folder')
# Unzip the file 
from zipfile import ZipFile
String='Ex_4_TR_PIV_Jet.zip'
zf = ZipFile(String,'r'); zf.extractall('./'); zf.close()
	\end{lstlisting}
\end{centering}

The dataset consists of $n_t=2000$ velocity fields, sampled at $f_s=2kHz$ over a grid of $60\times114$ points. The spatial resolution is approximately $\Delta x=0.33 mm$. The script provided in the book's website\footnote{\url{https://www.datadrivenfluidmechanics.com/chapter8}} describe how to manipulate this dataset, plot a time step and more. For the mPOD of this test case, the frequency splitting vector used in \cite{mendez2019multi} is constructed in terms of Strouhal number $St=f H/U_J$, that is the dimensionless frequency computed from the advection time $H/U_J$, with $H$ the stand-off distance and $U_J$ the mean velocity of the jet at the nozzle outlet. These are $H=4\mbox{mm}$ and $U_{J}=6.5m/s$. 

\begin{figure}[htbp]
	\centering
	\begin{minipage}{.49\linewidth}
		\includegraphics[keepaspectratio=true,width=1 \columnwidth]
		{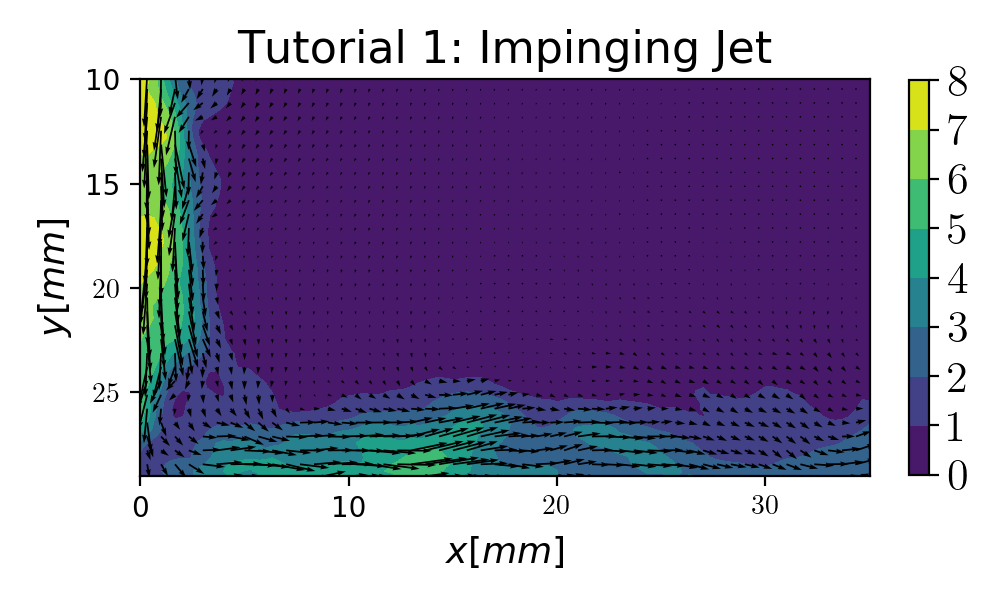}	
	\end{minipage}
	\begin{minipage}{.49\linewidth}
		\includegraphics[keepaspectratio=true,width=0.9 \columnwidth]
		{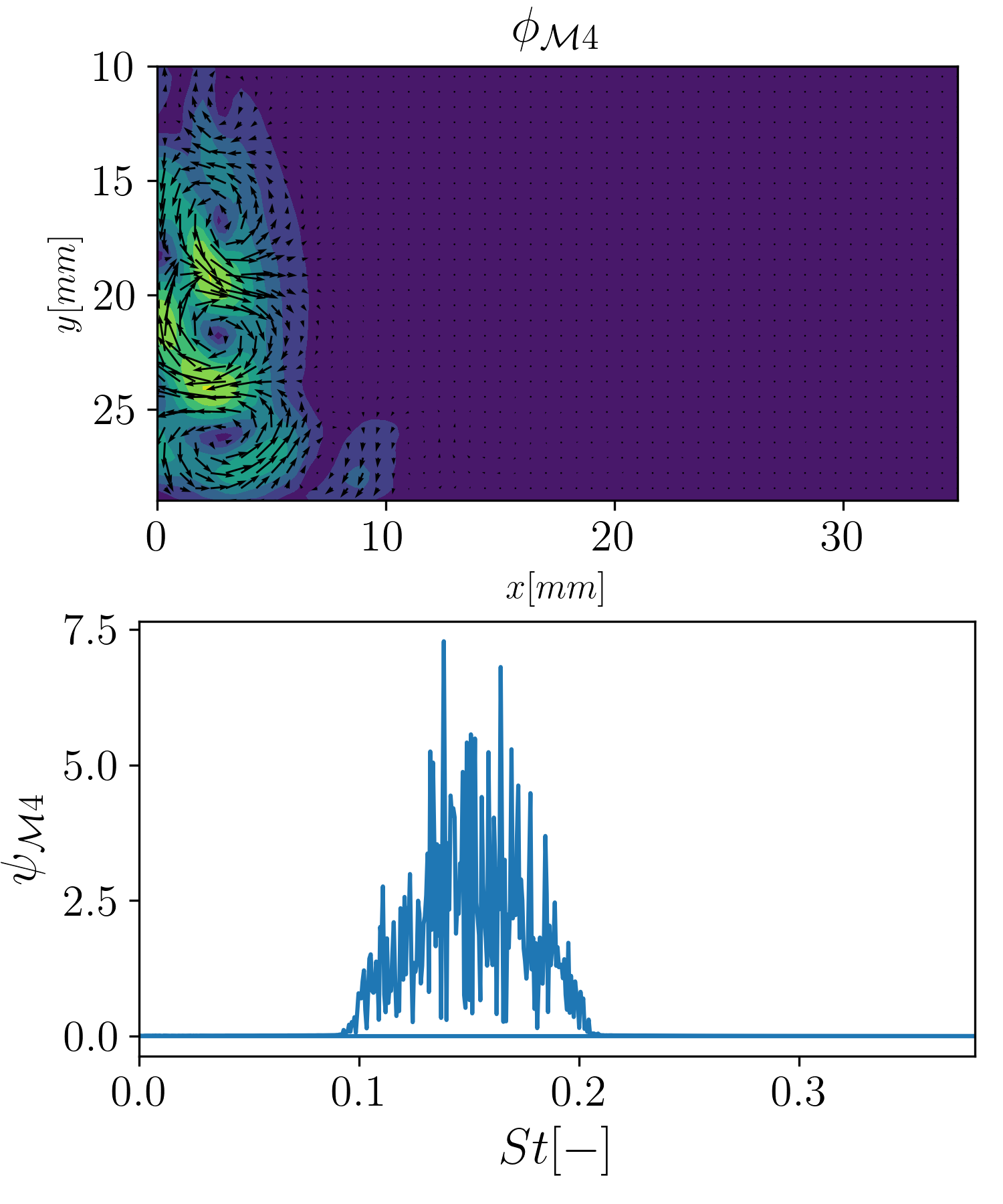}	
	\end{minipage}\\
	\caption{On the left: Example snapshot from the first tutorial test case on the TR-PIV of an impinging gas jet. Colormap in $m/s$. On the right: example of spatial structure (top) and frequency content (bottom) of the temporal structures of a mPOD mode.}\label{ch8fig13}
\end{figure}

The reader is encouraged to compare the mPOD results with those achievable from other decompositions. The spatial structures and the spectra of the temporal structures in the first $5$ modes are produced by the script \textsc{Tut\_1.py}. These are shown in Figure \ref{ch8fig13} for the dominant mPOD mode in the scale $St=0.1-0.2$. This mode isolates the large vortices roll-like structures produced as the instability of the shear layer evolve downstream the potential core of the jet. Recall that the velocity components are stacked into a single vector.

\subsection{Test Case 2}

To download this dataset and store it into a local directory, run the following script

\begin{centering}
	\begin{lstlisting}[language=Python,linewidth=12cm,xleftmargin=0\textwidth,xrightmargin=0\textwidth,backgroundcolor=\color{yellow!10}]
import urllib.request
print('Downloading Data for Tutorial 2...')
url = 'https://osf.io/qa8d4/download'
urllib.request.urlretrieve(url, 'Ex_5_TR_PIV_Cylinder.zip')
print('Download Completed! I prepare data Folder')
# Unzip the file 
from zipfile import ZipFile
String='Ex_5_TR_PIV_Cylinder.zip'
zf = ZipFile(String,'r'); 
zf.extractall('./DATA_CYLINDER'); zf.close()
	\end{lstlisting}
\end{centering}

\begin{figure}[htbp]
	\centering
	\begin{minipage}{.495\linewidth}
		\centering
		\includegraphics[keepaspectratio=true,width=1 \columnwidth]
		{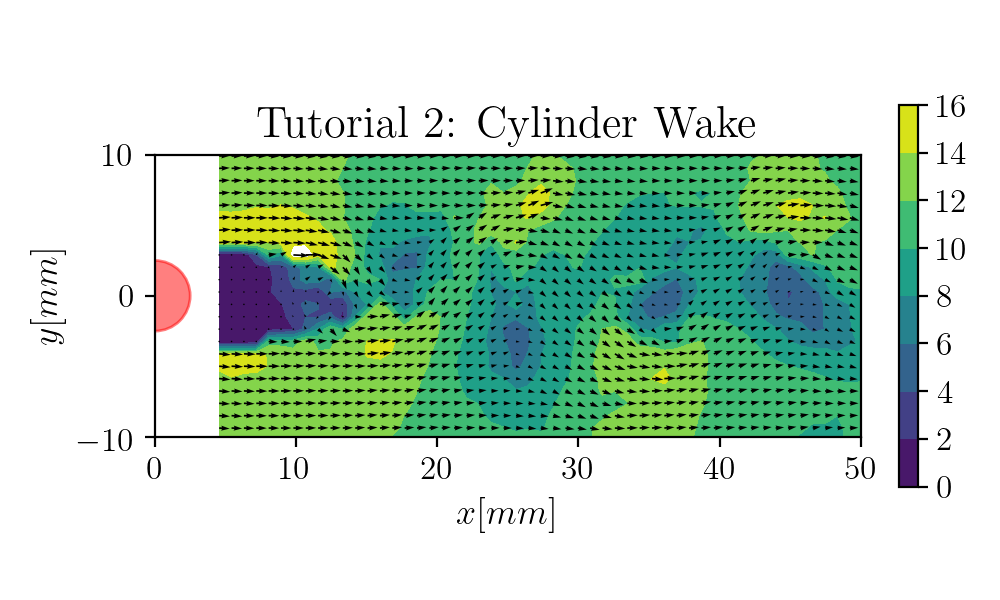}\\
		\vspace{-0.4cm}
		\includegraphics[keepaspectratio=true,width=1 \columnwidth]
		{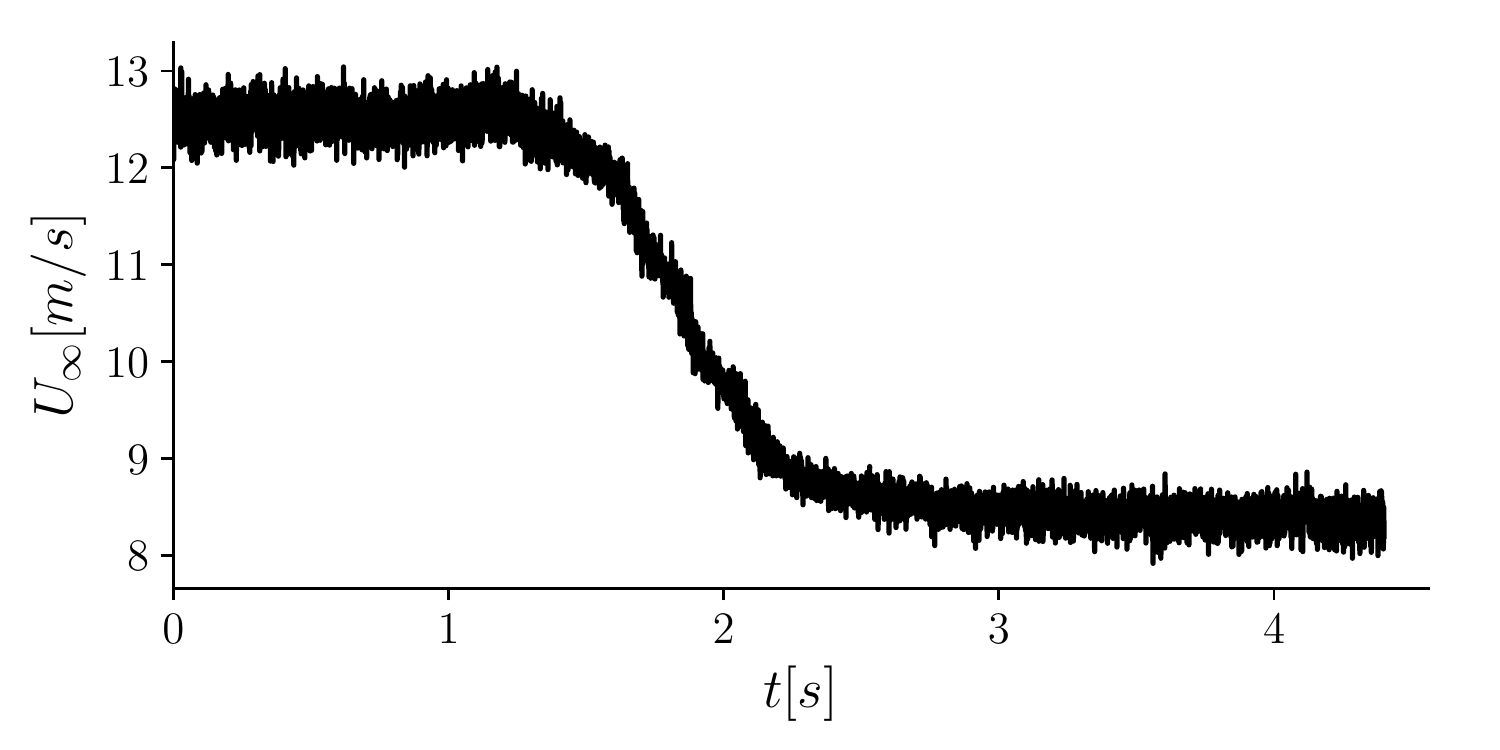}		
	\end{minipage}
	\begin{minipage}{.49\linewidth}
		\includegraphics[keepaspectratio=true,width=1 \columnwidth]
		{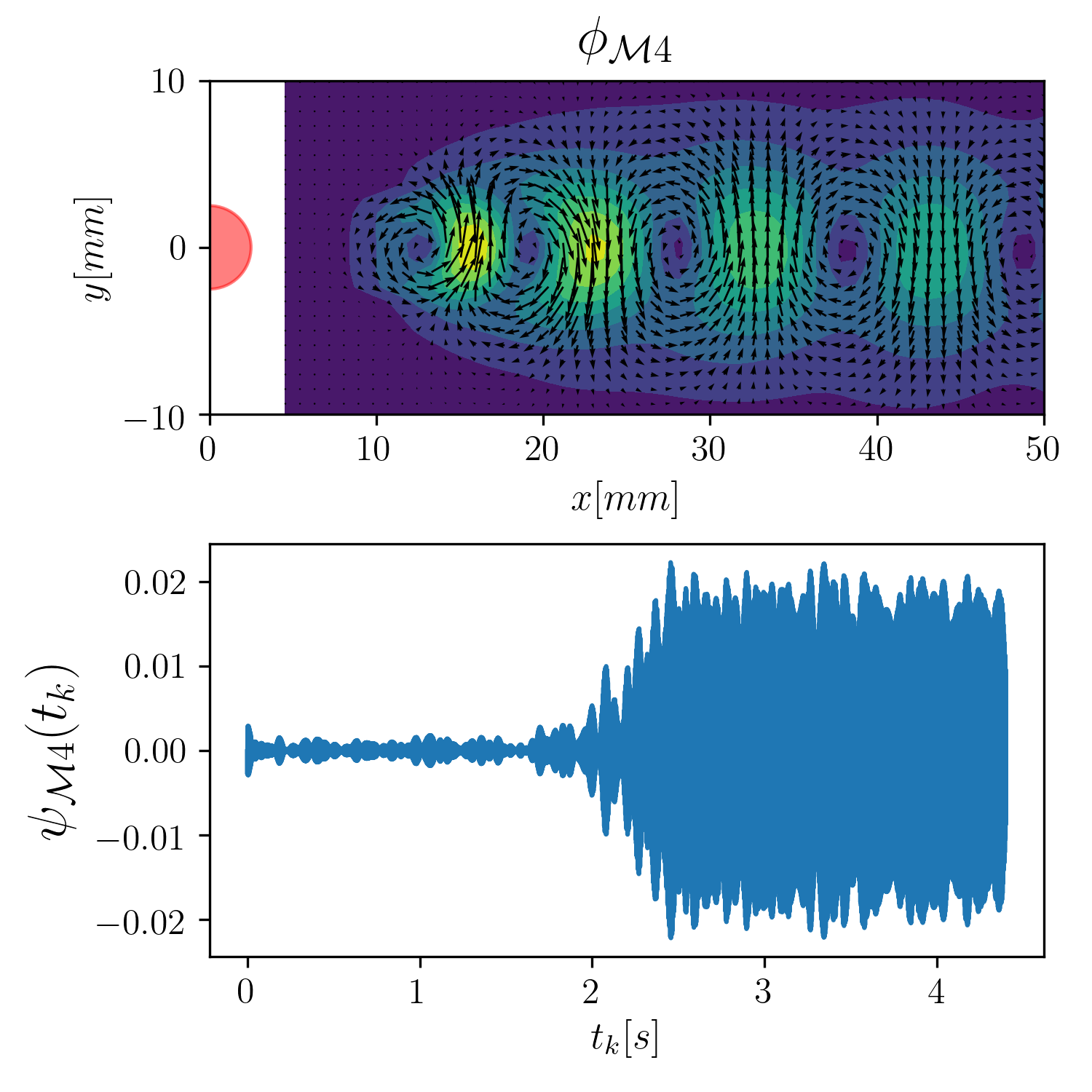}	
	\end{minipage}\\
	\vspace{-0.3cm}
	\caption
	{On the top left: a snapshot of the velocity field past a cylinder, obtained via TR-PIV. Colormap in $\mbox{m/s}$. On the bottom left: evolution of the free stream velocity as a function of time, sampled from the top left corner of the field. On the right: example of spatial structure (top) and frequency content (bottom) of the temporal structures of a mPOD mode.}\label{ch8fig14}
\end{figure}

The experimental set up and configuration is extensively described in \cite{DANTEC_VKI}. This dataset is computationally much more demanding than the previous, and consists of $n_t=13200$ velocity fields sampled at $f_s=3kHz$ over a grid of $71\times30$ points. The spatial resolution is approximately $\Delta x=0.85 mm$. A plot of a time step is shown in Figure \ref{ch8fig14}

As for the previous tutorial, a function is used to extract all the information about the grid (in this case stored on a different file) and the velocity field. 

This test case is characterized by a large scale variation of the free-stream velocity. A plot of the velocity magnitude in the main stream is shown in Figure \ref{ch8fig14}. In the first $1.5s$, the free stream velocity is at approximately $U_\infty\approx12\mbox{m/s}$. Between $t=1.5\mbox{s}$ and $t=2.5\mbox{s}$ this drops down to $U_\infty\approx8\mbox{m/s}$. The variation of the flow velocity is sufficiently low to let the vortex shedding adapt, and hence preserve an approximately constant Strhoual number of $St=f d U_{\infty}\approx 0.19$, with $d=5mm$ the diameter of the cylinder. Consequently, the vortex shedding varies from $f\approx 459 Hz$ to $f\approx 303 Hz$. Interestingly, the POD assign the entire evolution to a single pairs of modes, and it is hence not possible to analyze the vortex structures in the shedding for the two phases at approximately constant velocity, nor to distinguish them from the flow organization during the transitory phase.

The mPOD can be used to identify modes related to these three phases.
Three scales are chosen in the exercise. The first, in the range $\Delta f=[0-10] \mbox{Hz}$ is designed to isolate the large scale motion of the flow, hence the variation of the free stream velocity. The scale with $\Delta f=[290-320]\mbox{Hz}$ is designed to capture the vortex shedding when the velocity is at $U_\infty\approx8\mbox{m/s}$ while the scale with $\Delta f=[430-470]\mbox{Hz}$ identifies the shedding when the velocity is at $U_\infty\approx12\mbox{m/s}$.

The spatial structure and the temporal evolution of am mPOD mode is also shown in Figure \ref{ch8fig14}. This mode is clearly associated to the second stationary conditions. Other modes captures the first part while others are focused on the transitory phase: besides allowing for spectral localization, the MRA structure of the decomposition also allows for providing time localization capabilities.

\section{What's next?}\label{sec7}

This chapter presented a general matrix factorization framework to formulate \emph{any} linear decomposition. This general framework allowed for deriving the Multiscale Proper Orthogonal Decomposition, a formulation that allows to bridge energy-based and frequency-based formulations. The proposed exercises and the extensive collection of codes listed in the book's website should allow the reader to efficiently use these decompositions on his/her own data and provide a good starting point to develop computational proficiency on this subject.

It is important to note that this chapter has only focused on the mathematical framework of \emph{linear decompositions}: every decomposition provides its own way of representing the data with respect to a certain basis. From here, the reader might be interested in two possible directions. On the representation side, while linear dimensionality reduction has been at the heart of an extensive body of literature in fluid mechanics, nonlinear methods are now offering a new set of tools as show-cased in Chapter 1. On the reduced order modeling side, whether a low dimensional representation is capable of reducing the complexity of a dynamical system depends on whether the projection (or the nonlinear mapping) is able to preserve the most dynamically relevant phenomena in the new representation. These challenges are described in Chapter 14.

\bibliographystyle{apalike}
\bibliography{Chapter_8_DDFM_Book_Mendez}
\end{document}